\documentstyle[amssymb]{article}

\newtheorem{theorem}{Theorem}

\newtheorem{definition}[theorem]{Definition}
\newtheorem{example}[theorem]{Example}

\def\QED{\quad\blackslug\lower 8.5pt\null}

\begin{document}

\begin{center}
{\Large \bf THE GEOMETRY } 
\vspace*{2mm}

{\Large \bf OF LIGHTLIKE  HYPERSURFACES}

\vspace*{2mm}

{\Large \bf OF THE DE SITTER SPACE}

\vspace*{3mm}
{\large M.A. Akivis and  V.V. Goldberg}
\end{center}

\vspace*{5mm}

{\footnotesize 
{\it Abstract}:
It is proved that the geometry of lightlike 
hypersurfaces of the de Sitter space 
$S^{n+1}_1$ is directly connected with the 
geometry of hypersurfaces of the conformal space 
$C^n$. This connection is applied for 
a construction of an invariant normalization and an
invariant affine connection  of lightlike 
hypersurfaces as well as for studying singularities 
of lightlike hypersurfaces.

\vspace*{2mm}


distribution, singularity, 
affine connection.

\vspace*{2mm}

{\it MS Classification}: 53B30, 53A30,  53B50, 53A35, 
\setcounter{equation}{0}
}
\setcounter{section}{-1}

\section{Introduction}  

The projective model of the non-Euclidean geometry 
(the Cayley--Klein model) is closely connected with models of 
conformal geometry and the geometry of the de Sitter space. 
In fact, the hyperbolic space $H^{n+1}$ of dimension $n+1$---
the Lobachevsky space---admits a mapping onto internal 
domain of an $n$-dimensional oval hyperquadric $Q^n$ of 
a projective space $P^{n+1}$. On this hyperquadric itself 
the geometry of an $n$-dimensional conformal space $C^n$ 
is realized, and outside of the hyperquadric $Q^n$ 
the geometry of the $(n+1)$-dimensional  de Sitter space 
$S^{n+1}_1$ 
is realized. Moreover, the group of projective transformations 
of the  space $P^{n+1}$ keeping the hyperquadric $Q^n$ invariant 
and transferring its internal domain into itself 
(this group is denoted by ${\bf PO} (n+2, 1)$---see 
\cite{AG96}, p. 7) 
is isomorphic to the group of motions of the Lobachevsky space 
$H^{n+1}$, the conformal space $C^n$, and the de Sitter space 
$S^{n+1}_1$. It is clear that there exist deep connections 
among these three geometries.

The Lobachevsky geometry is the first example of geometry 
which differs from the Euclidean geometry. Numerous books 
and papers are devoted to the Lobachevsky geometry. 
Conformal differential geometry was also studied in detail. 
In particular, it was studied in the last authors' book \cite{AG96}. In spite of the fact that the geometry of the de 
Sitter space is the simplest model of spacetime of 
general relativity, this geometry was not studied thoroughly. 
The de Sitter space sustains the Lorentzian metric of constant 
positive curvature. 

In the present paper we study the geometry of the 
de Sitter space $S^{n+1}_1$ using its connection with the 
geometry of the conformal space. We prove that the geometry of 
lightlike hypersurfaces of the space $S^{n+1}_1$, 
which play an important role in general relativity 
(see the book {\cite{DB96}), is directly connected with the 
geometry of hypersurfaces of the conformal space $C^n$. The 
latter  was studied in detail in the papers of the first 
author (see \cite{A52}, \cite{A61}, \cite{A63}, \cite{A64}, 
\cite{A65}) and also in the book \cite{AG96}. This simplifies the 
study of lightlike hypersurfaces of the de Sitter space 
$S^{n+1}_1$ and makes possible to apply for their consideration 
the apparatus constructed in the conformal theory. 

In Section 1 we study the geometry of the 
de Sitter space and its connection with the geometry of the 
conformal space. After this we study lightlike hypersurfaces 
$U^n$ in the space $S^{n+1}_1$, investigate their structure, 
and prove that such a hypersurface is tangentially degenerate of 
rank $r \leq n - 1$. Its rectilinear or plane 
generators form an isotropic fibre bundle on $U^n$.

In Sections 2--5 we investigate lightlike hypersurfaces $U^n$ of 
maximal rank, and for their study we use the relationship between 
the geometry of such hypersurfaces and the geometry of 
hypersurfaces of the conformal space. For a lightlike 
hypersurface, we construct the fundamental quadratic forms and 
connections determined  by a normalization of a hypersurface by 
means of a distribution (the screen distribution) which is 
complementary to the isotropic distribution. The screen 
distribution plays an important role in the book \cite{DB96} 
since it defines a connection on a lightlike hypersurface $U^n$, 
and it appears to be important for applications. We prove that 
the screen distribution on a lightlike hypersurface  can be 
constructed invariantly by means of quantities from 
a third-order differential neighborhood, that is, such a distribution is 
intrinsically connected with the geometry of a hypersurface. 

In Section 5 we study singular points of a  lightlike 
hypersurface in the de Sitter space $S^{n+1}_1$, classify them, 
and describe the structure of hypersurfaces carrying singular 
points of different types. Moreover, we establish the connection 
of this classification with that of canal hypersurfaces of the 
conformal space.

In Section 6 we consider lightlike hypersurfaces of reduced rank. 
Such hypersurfaces carry lightlike rectilinear generators along 
which their tangent hyperplanes are constant. For such 
hypersurfaces, again in a third-order differential neighborhood 
we construct an invariant screen distribution and an invariant 
affine connection. However, the method of construction 
is different from that for lightlike hypersurfaces of maximal 
rank, 
since the construction used for hypersurfaces of maximal rank 
fails for hypersurfaces of reduced rank. We establish a 
connection of lightlike hypersurfaces of reduced rank 
with quadratic hyperbands of 
a multidimensional projective space.

The principal method of our investigation is the method of 
moving frames and exterior differential forms in the form 
in which it is presented in the books \cite{AG93} and 
\cite{AG96}. All functions considered in the paper are assumed to 
be real and 
differentiable, and all manifolds are assumed to be smooth 
with the possible exception of some isolated singular points and 
singular submanifolds.

\section{The de Sitter Space}

\setcounter{equation}{0}
 
{\bf 1.} 
In a projective space $P^{n+1}$ of dimension $n+1$ we consider 
an oval hyperquadric $Q^n$. Let $x$ be a point of the space 
$P^{n+1}$ with projective coordinates $(x^0, x^1, \ldots , 
x^{n+1})$.  The hyperquadric $Q^n$ is determined by the equations 
\begin{equation}\label{eq:1}
(x, x) := g_{\xi\eta} x^\xi x^\eta = 0, \;\;\;\;\; 
\xi, \eta = 0, \ldots, n+1,
\end{equation}
whose left-hand side is a quadratic form $(x, x)$ 
of signature $(n+1, 1)$. The hyperquadric $Q^n$ 
divides the space $P^{n+1}$ into two parts, external and 
internal. Normalize the quadratic form $(x, x)$ 
in such a way that for the points of the external part 
the inequality $(x, x) > 0$ holds. This external domain is 
a model of the {\em de Sitter space} $S^{n+1}_1$ 
(see \cite{Z96}). We will identify the external domain 
of $Q^n$ with the space $S^{n+1}_1$. The hyperquadric $Q^n$ 
is the {\em absolute} of the space $S^{n+1}_1$.

On the hyperquadric $Q^n$ of  the space $P^{n+1}$ the geometry 
of a conformal space $C^n$ is realized. The bijective mapping 
$C^n \leftrightarrow Q^n$ 
is called the {\em Darboux mapping}, and 
the hyperquadric $Q^n$ itself is called the {\em Darboux 
hyperquadric}. 

Under the Darboux mapping to hyperspheres of the space $C^n$ 
there correspond cross-sections of the hyperquadric $Q^n$ 
by hyperplanes $\xi$. But to a hyperplane $\xi$ there corresponds 
a point $x$ that is polar-conjugate to $\xi$ with 
respect to $Q^n$ and lies outside of $Q^n$, that is, 
a point of the space $S^{n+1}_1$. Thus to hyperspheres of the 
space $C^n$ there correspond points  of the space $S^{n+1}_1$. 

Let $x$ be an arbitrary point of the space $S^{n+1}_1$. 
The tangent lines from the point $x$ to the hyperquadric $Q^n$ 
form a second-order cone $C_x$ with vertex at the point $x$. 
This cone is called the {\em isotropic cone}. 
For spacetime whose model 
is the space $S^{n+1}_1$ this cone is the light cone, and 
its generators are lines of propagation of light impulses 
whose source coincides with  the point $x$. 

The cone $C_x$ separates all straight lines passing through 
 the point $x$ into spacelike (not having common points with 
the hyperquadric $Q^n$), timelike (intersecting $Q^n$ 
in two different points), and lightlike (tangent to $Q^n$). 
The  lightlike straight lines are generators of the cone $C_x$. 

To a spacelike straight line $l \subset S_1^{n+1}$ there 
corresponds an elliptic pencil of hyperspheres in the conformal 
space $C^n$. All hyperspheres of this pencil pass through a 
common $(n-2)$-sphere $S^{n-2}$ (the center of this pencil). 
The sphere $S^{n-2}$ is the intersection of the hyperquadric 
$Q^n$ and the $(n-1)$-dimensional subspace of the space $P^{n+1}$ 
which is polar-conjugate to the line $l$ with respect to the 
hyperquadric $Q^n$. 

To a timelike straight line $l \subset 
S_1^{n+1}$ there corresponds a hyperbolic pencil of hyperspheres 
in the space $C^n$. 
Two arbitrary hyperspheres of this pencil do 
not have common points, and the pencil contains two hyperspheres 
of zero radius which correspond to the points of intersection of 
the straight line $l$ and the hyperquadric $Q^n$. 

Finally, to a lightlike straight line $l 
\subset S_1^{n+1}$ there corresponds a parabolic pencil of 
hyperspheres in the  space $C^n$ consisting of 
hyperspheres tangent one to another at a point that is a unique 
hypersphere of zero radius belonging to this pencil.

Hyperplanes of the space $S^{n+1}_1$ are also divided 
into three types. Spacelike hyperplanes do not have 
common points with  the hyperquadric $Q^n$; a timelike 
hyperplane intersects   $Q^n$ along a real hypersphere; 
and lightlike hyperplanes are tangent to $Q^n$. 
Subspaces of any dimension $r, \; 2 \leq r \leq n-1$, 
can be also classified in a similar manner.

Let us apply the method of moving frames 
to study some questions of differential geometry 
of the space $S^{n+1}_1$. With a point $x \in S^{n+1}_1$ 
we associate a family of projective frames $\{A_0, A_1, 
\ldots , A_{n+1}\}$. However, in order to apply formulas derived 
in the book \cite{AG96}, we will use the notations used 
in this book. Namely, we denote by $A_n$ the vertex 
of the moving frame which coincides 
 with the point $x$, $A_n = x$; we locate the vertices 
$A_0, A_i$, and $A_{n+1}$ at the hyperplane $\xi$ which is polar 
conjugate to the point $x$ with respect to 
 the hyperquadric $Q^n$, and we assume that the points $A_0$ and 
$A_{n+1}$ lie on the hypersphere $S^{n-1} = Q^n \cap \xi$, and 
the points $A_i$ are polar-conjugate to the straight line $A_0 
A_{n+1}$ with respect to $S^{n-1}$. Since $(x, x) > 0$, we can 
normalize the point $A_n$ by the condition $(A_n, A_n) = 1$. 
The points  $A_0$ and $A_{n+1}$ are not polar-conjugate  with 
respect to  the hyperquadric $Q^n$. Hence we can normalize them 
by the condition $(A_0, A_{n+1}) = - 1$. As a result, 
the matrix of scalar products of the frame elements 
has the form 
\begin{equation}\label{eq:2}
\renewcommand{\arraystretch}{1.3}
(A_\xi, A_\eta) = \pmatrix{0 & 0 & 0 & -1 \cr
0 & g_{ij} & 0 & 0 \cr
 0 & 0 & 1 & 0 \cr
-1 & 0 & 0& 0}, \;\;\;\;\; i, j = 1, \ldots, n -1,
\renewcommand{\arraystretch}{1}
\end{equation}
and the quadratic form $(x, x)$ takes the form
\begin{equation}\label{eq:3}
(x, x) = g_{ij} x^i x^j + (x^n)^2 - 2 x^0 x^{n+1}.
\end{equation}
 The quadratic form $g_{ij} x^i x^j $ occurring in 
(3) is positive definite. 

The equations of infinitesimal displacement of the conformal 
frame $\{A_\xi\}, \linebreak \xi = 0, 1, \ldots, n+1$, 
we have constructed have the form
\begin{equation}\label{eq:4}
d A_\xi = \omega_\xi^\eta A_\eta, \;\;\;\; \xi, \eta 
= 0, 1, \ldots, n+1,
\end{equation}
where by (2), the 1-forms $\omega_\xi^\eta$  
  satisfy the following Pfaffian equations: 
\begin{equation}\label{eq:5}
\left\{
\renewcommand{\arraystretch}{1.3}
\begin{array}{ll}
\omega_0^{n+1} =  \omega^0_{n+1} = 0, &  \omega_0^0  + \omega_{n+1}^{n+1} = 0, \\ 
 \omega_i^{n+1} = g_{ij} \omega_0^j, & 
\omega_i^0 = g_{ij} \omega_{n+1}^j, \\ 
\omega_n^{n+1} - \omega_0^n   = 0, & \omega_n^0 - \omega_{n+1}^n = 0, \\
 g_{ij} \omega_n^j + \omega_i^n = 0 , &   \omega_n^n = 0, \\ 
 dg_{ij} = g_{jk} \omega_i^k + g_{ik} \omega_j^k. 
\end{array}
\renewcommand{\arraystretch}{1}
\right.
\end{equation}
These formulas are precisely the formulas derived in the book 
\cite{AG96} (see p. 32) for the conformal space $C^n$. 

 It follows from (4) that 
\begin{equation}\label{eq:6}
dA_n = \omega_n^0 A_0 + \omega_n^i A_i + \omega_n^{n+1} A_{n+1}.
\end{equation}

The differential $dA_n$ belong to the tangent space 
$T_x (S_1^{n+1})$, and the 1-forms $\omega_n^0, \omega_n^i$, 
and $\omega_n^{n+1}$ form a coframe of this space. The total 
number of these forms is $n + 1$, and this number coincides with 
the dimension of $T_x (S_1^{n+1})$. The scalar square of the 
differential $dA_n$ is the metric quadratic form $\widetilde{g}$ on 
the manifold $S_1^{n+1}$. By (2), this quadratic form $\widetilde{g}$ 
can be written as 
$$
\widetilde{g} = (dA_n, dA_n) = g_{ij} \omega_n^i \omega_n^j - 2 \omega_n^0 
\omega_n^{n+1}.
$$
Since the first term of this expression is a positive definite 
quadratic form, the form $\widetilde{g}$ is of Lorentzian signature 
$(n, 1)$. The coefficients of the form $\widetilde{g}$ produce the metric tensor 
of the space $S^{n+1}_1$ whose matrix is obtained from the 
matrix (2) by deleting the $n$th row and the $n$th column. 

The quadratic form $\widetilde{g}$  defines 
on $S^{n+1}_1$ a pseudo-Riemannian metric of signature $(n, 1)$. 
The isotropic cone defined in the space 
$T_x (S_1^{n+1})$ by the equation $\widetilde{g} = 0$ coincides with the 
cone $C_x$ that we defined earlier in the space $S^{n+1}_1$ 
geometrically. 

The 1-forms $\omega_\xi^\eta$ occurring in equations (4) 
  satisfy the structure equations of the space $C^n$:
\begin{equation}\label{eq:7}
d \omega_\xi^\eta = \omega_\xi^\zeta \wedge \omega_\zeta^\eta,
\end{equation}
which are obtained by taking exterior derivatives of 
equations (4) and which are conditions of complete 
integrability of (4). The forms $\omega_\xi^\eta$ 
are invariant forms of the fundamental group 
${\bf PO} (n+2, 1)$ of transformations of the 
spaces $H^{n+1}, C^n$, and $S^{n+1}_1$ which is locally 
isomorphic to the group ${\bf SO} (n+2, 1)$. 

Let us write equations (7) for the 1-forms 
$\omega_n^0, \omega_n^i$, and $\omega_n^{n+1}$ 
making up a coframe of the space $T_x (S_1^{n+1})$  
in more detail: 
\begin{equation}\label{eq:8}
\renewcommand{\arraystretch}{1.3}
\begin{array}{lll}
d \omega_n^0 = & \omega_n^0 \wedge \omega_0^0 
                  + &\omega_n^i \wedge \omega_i^0, \\
d \omega_n^i = & \omega_n^0 \wedge \omega_0^i 
                  + &\omega_n^j \wedge \omega_j^i 
                  + \omega_n^{n+1} \wedge \omega_{n+1}^i, \\
d \omega_n^{n+1} = & &\omega_n^i \wedge \omega_i^{n+1} 
                  + \omega_n^{n+1} \wedge \omega_{n+1}^{n+1}.
\end{array}
\renewcommand{\arraystretch}{1}
\end{equation}
The last equations can be written in the matrix form 
as follows:
\begin{equation}\label{eq:9}
d \theta = - \omega \wedge \theta,
\end{equation}
where $\theta = (\omega_n^u), u = 0, i, n+1$, is 
the column matrix with its values in the vector space 
 $T_x (S_1^{n+1})$, and $\omega = (\omega_v^u), u, v = 0, i, n+1$, is a square matrix of order $n + 1$ with values in 
the Lie algebra of the group of admissible transformations of 
coframes of the space  $T_x (S_1^{n+1})$. The form $\omega$ 
is the connection form of the space  $S_1^{n+1}$. 
In detail this form can be written  as
\begin{equation}\label{eq:10}
 \omega =\pmatrix{ \omega_0^0 &  \omega_i^0 & 0 \cr 
                  \omega_0^i &  \omega_i^j &  \omega_{n+1}^i \cr
               0 &  \omega_i^{n+1} &  \omega_{n+1}^{n+1}}.
\end{equation}
By (5), in this matrix, only the forms in the left 
upper corner, which form an $n\times n$-matrix, 
are linearly independent.

The connection form (10) allows us to find the differential 
equation of geodesics in the space $S^{n+1}_1$. These lines 
coincide with straight lines of the ambient space $P^{n+1}$; 
more precisely, they coincide with the parts of these straight 
lines which lie outside of the Darboux hyperquadric $Q^n$. 
We will look for their equation in the form $x = x (t)$, 
and we will impose the vertex $A_n$ of the moving frame with 
the point $x$, $A_n = x (t)$. Write the decomposition 
of the tangent vector to a geodesic in the form
$$
\frac{dx}{dt} = \xi^u A_u, \;\;\;\;\; u = 0, i, n+1.
$$
For a geodesic, the second differential 
$\frac{d^2x}{dt^2}$ is collinear to its tangent vector 
$\frac{dx}{dt}$. This implies that 
$$
\frac{d\xi^u}{dt} A_u + \xi^v \omega_v^u A_u = \alpha \xi^u A_u,
$$
where the connection 1-forms $\omega_v^u$ composing the matrix 
(10) are calculated along the curve $x = x (t)$, and $\alpha$ is 
a new 1-form. Hence the differential equation of geodesics 
has the form
\begin{equation}\label{eq:11}
\frac{d\xi^u}{dt}  + \xi^v \omega_v^u  = \alpha \xi^u.
\end{equation}
The same equation (11) is the equation of straight lines 
of the space $P^{n+1}$.

Next we will find the curvature form and  the curvature 
tensor of the space $S^{n+1}_1$. To this end, we 
take exterior derivative of the connection form $\omega$, 
more precisely, of its independent part. Applying 
equations (7), we find the following components of the 
 curvature form: 
\begin{equation}\label{eq:12}
\renewcommand{\arraystretch}{1.3}
\left\{
\begin{array}{lll}
\Omega_0^0 = d \omega_0^0 - \omega_0^i \wedge \omega_i^0 
                  = \omega_n^{n+1} \wedge \omega_n^0, \\
\Omega_0^i = d \omega_0^i - \omega_0^0 \wedge \omega_0^i 
                   -  \omega_0^j \wedge \omega_j^i 
=   \omega_n^{n+1} \wedge \omega_n^i, \\
\Omega_i^0 = d \omega_i^0 - \omega_i^0 \wedge \omega_0^0 
                   -  \omega_i^j \wedge \omega_j^0 
= -g_{ij}  \omega_n^j \wedge \omega_n^0, \\
\Omega_j^i = d \omega_j^i - \omega_j^0 \wedge \omega_0^i 
                   -  \omega_j^k \wedge \omega_k^i  
                   -  \omega_j^{n+1} \wedge \omega_{n+1}^i  
= -g_{jk}  \omega_n^k \wedge \omega_n^i. 
\end{array}
\renewcommand{\arraystretch}{1}
\right.
\end{equation}

But the general expression of the  curvature form of 
an $(n+1)$-dimensional pseudo-Riemannian space with 
a coframe $\omega_n^0, \omega_n^i,$ and $\omega_n^{n+1}$ 
has the form 
\begin{equation}\label{eq:13}
 \Omega_s^r = d\omega_s^r -  \omega_s^t \wedge 
               \omega_t^r = 
\frac{1}{2} R^r_{suv}  \omega_n^u \wedge  \omega_n^v, 
\end{equation}
where $r, s, t, u, v = 0, 1, 
\ldots , n -1, n + 1$ (see, for example, \cite{W77}). 
Comparing equations (12) and (13), we find that 
$$
\Omega_s^r = \omega_u^n \wedge g_{sv} \omega_n^v
$$
and
\begin{equation}\label{eq:14}
R^r_{suv} = \delta_u^r g_{sv} - \delta_v^r g_{su},
\end{equation}
where $(g_{sv})$ is the matrix of coefficients of the quadratic form (2). But this means that the space $S^{n+1}_1$ 
is a pseudo-Riemannian space of constant positive curvature $K = 1$. The Ricci tensor of this space has the form 
\begin{equation}\label{eq:15}
R_{sv} = R^r_{srv} = n  g_{sv}.
\end{equation}
This confirms that the space $S^{n+1}_1$, as any 
pseudo-Riemannian space of constant  curvature, 
is the Einstein space.

Thus by means of the method of moving frame we proved the following well-known theorem (see, for example, \cite{W77}): 

\begin{theorem} The de Sitter space, whose model is 
the domain of a projective space $P^{n+1}$ lying outside 
of an oval hyperquadric $Q^n$, is a pseudo-Riemannian space of 
Lorentzian signature $(n, 1)$ and of constant positive curvature 
$K = 1$. This space is homogeneous, and its fundamental group 
${\bf PO} (n+2, 1)$ is locally isomorphic 
to the special orthogonal group ${\bf SO} (n+2, 1)$. 
 \end{theorem}

\section{Lightlike Hypersurfaces in the de Sitter Space}

A hypersurface $U^n$ in the de Sitter space $S^{n+1}_1$ 
is said to be {\em lightlike} if all its tangent hyperplanes are 
lightlike, that is, they are tangent to the hyperquadric $Q^n$ 
which is the absolute of the  space $S^{n+1}_1$.

Denote by $x$ an arbitrary point of the  hypersurface $U^n$, 
by $\eta$ the tangent hyperplane to $U^n$ at the point $x,  
\eta = T_x (U^n)$, and by $y$ the point of tangency of 
the hyperplane $\eta$ with the hyperquadric $Q^n$. 
Next, as in Section 1,  denote by $\xi$ the hyperplane 
which is polar-conjugate to the point $x$ 
with respect to the hyperquadric $Q^n$, and  associate 
with a point $x$ a family of projective  frames such that 
$x = A_n, y = A_0$, the points $A_i, i = 1, \ldots , 
n - 1$, belong to the intersection of the hyperplanes $\xi$ and 
$\eta, \; A_i \in \xi \cap \eta$, and the point $A_{n+1}$, as 
well as the point $A_0$, belong to the straight line that is polar-conjugate to the $(n-2)$-dimensional subspace spanned by the points $A_i$. In addition, we normalize the frame vertices 
in the same way as this was done in Section 1. Then the matrix of scalar products of the frame elements has the form (2), and the components of infinitesimal displacements of the moving frame 
satisfy the Pfaffian equations (5). 

Since the hyperplane $\eta$ is tangent to 
the  hypersurface $U^n$ at the point $x = A_n$ and does not 
contain the point $A_{n+1}$, the differential of 
 the point $x = A_n$ has the form
\begin{equation}\label{eq:16}
dA_n =\omega_n^0 A_0 + \omega_n^i A_i,
\end{equation}
 the following equation holds:
\begin{equation}\label{eq:17}
\omega_n^{n+1} = 0,
\end{equation}
and the forms $\omega_n^0$ and $\omega_n^i$ are basis forms of 
the  hypersurface $U^n$.
By relations (5), it follows from equation (16) that 
\begin{equation}\label{eq:18}
\omega_0^n = 0
\end{equation}
and 
\begin{equation}\label{eq:19}
dA_0 =\omega_0^0 A_0 + \omega_0^i A_i.
\end{equation}

Taking  exterior derivative of equation (17), we obtain 
$$
\omega^i_n \wedge \omega_i^{n+1} = 0.
$$
Since the forms $\omega^i_n$ are linearly independent, 
by Cartan's lemma, we find from the last equation that 
\begin{equation}\label{eq:20}
\omega_i^{n+1} = \nu_{ij} \omega_n^j, \;\; 
\nu_{ij} = \nu_{ji}.
\end{equation}
Applying an appropriate formula from (5), we find that
\begin{equation}\label{eq:21}
\omega_0^i = g^{ij}  \omega_j^{n+1} = 
g^{ik} \nu_{kj} \omega_n^j,
\end{equation}
where $(g^{ij})$ is the inverse matrix of the matrix $(g_{ij})$. 

Now formulas (16) and (19) imply that for 
$\omega_n^i = 0$, the point $A_n$ of 
the  hypersurface $U^n$ moves along the isotropic straight line 
$A_n A_0$, and hence $U^n$ is a ruled hypersurface. In what 
follows, we assume that the {\em entire}  straight line 
$A_n A_0$ belongs to the  hypersurface $U^n$. 

In addition, formulas (16) and (19) show that at any point of 
a generator of  the  hypersurface $U^n$, its tangent hyperplane 
is fixed and coincides with the hyperplane $\eta$. Thus 
   $U^n$ is a {\em tangentially degenerate   hypersurface}. 

We recall that the {\em rank} of a tangentially degenerate   
hypersurface is the number of parameters on which the family of 
its tangent hyperplanes depends (see, for example, \cite{AG93}, 
p. 113). From relations (16) and (19) it follows that the 
tangent hyperplane $\eta$ of the hypersurface $U^n$ along 
its generator $A_n A_0$ is determined by this generator 
and the points $A_i$, $\eta = A_n \wedge A_0 \wedge A_1 
\wedge \ldots \wedge A_{n-1}$. The displacement of 
this hyperplane is determined by the differentials 
(16), (19), and
$$
dA_i = \omega_i^0 A_0 + \omega_i^j A_j + \omega_i^n A_n 
+ \omega_i^{n+1} A_{n+1}.
$$
But by (5), $\omega_i^n = - g_{ij} \omega_n^j$, and the forms 
$\omega_i^{n+1}$ are expressed according to formulas (20). 
From formulas (20) and (21) it follows that the rank of 
 a tangentially degenerate   hypersurface $U^n$ is 
determined by the rank of the matrix $(\nu_{ij})$ in 
terms of which the 1-forms $\omega_i^{n+1}$ and $\omega_0^i$ 
are expressed. But by (19) and (21) the dimension of 
the submanifold $V$ described by the point $A_0$ on the 
hyperquadric $Q^n$ is also equal to the rank of the matrix 
$(\nu_{ij})$. Thus we have proved 
the following result:

\begin{theorem} A lightlike hypersurface of the de Sitter 
space $S^{n+1}_1$ is a ruled tangentially degenerate 
hypersurface whose rank is equal to the dimension of 
the submanifold $V$ described by the point $A_0$ on the 
hyperquadric $Q^n$. 
\end{theorem}

Denote the rank of the tensor $\nu_{ij}$ and of  the hypersurface 
$U^n$ by $r$.  In this and next sections  we will assume that 
$r  = n - 1$, and 
the case $r < n -1$ will be considered in the last section of the paper.

For  $r = n - 1$, the hypersurface $U^n$ carries an 
$(n-1)$-parameter 
family of isotropic rectilinear generators $l = A_n A_0$ along 
which the tangent hyperplane $T_x (U^n)$ is fixed. From the point 
of view of physics, the isotropic rectilinear generators of 
a lightlike hypersurface $U^n$ are trajectories of light impulses, and 
the hypersurface $U^n$ itself represents a {\em light flux} in spacetime.

Since $\mbox{{\rm rank}} \; (\nu_{ij}) = n - 1$,
the submanifold $V$ described by the point $A_0$ on the 
hyperquadric $Q^n$ has dimension $n-1$, that is, 
$V$ is a hypersurface. We denote it by $V^{n-1}$. The tangent 
subspace $T_{A_0} (V^{n-1})$ to $V^{n-1}$ is determined by 
the points $A_0, A_1, \ldots , A_{n-1}$. Since $(A_n, A_i) 
= 0$, this tangent subspace is polar-conjugate to the 
rectilinear generator $A_0 A_n$ of 
the lightlike hypersurface $U^n$.

The submanifold $V^{n-1}$ of the 
hyperquadric $Q^n$ is the image of a hypersurface of 
the conformal space $C^n$ under the Darboux mapping. 
 We will denote this hypersurface also by $V^{n-1}$. 
In the space $C^n$, 
the hypersurface  $V^{n-1}$ is defined by equation (18) which 
by (5) is equivalent to equation (17) defining a lightlike 
hypersurface $U^n$ in the space $S_1^{n+1}$. To 
the rectilinear generator $A_n A_0$ of the hypersurface $U^n$ 
there corresponds a parabolic pencil of hyperspheres $A_n + sA_0$ 
tangent to the hypersurface  $V^{n-1}$ (see \cite{AG96}, p. 40). 
Thus the following theorem is valid:

\begin{theorem} There exists a one-to-one correspondence 
between the set of hypersurfaces of the conformal space $C^n$ 
and the set of lightlike hypersurfaces of the maximal rank 
$r = n-1$ of the de Sitter 
space $S^{n+1}_1$. To  pencils of tangent hyperspheres 
of the hypersurface  $V^{n-1}$ there correspond isotropic 
rectilinear generators of the lightlike hypersurface $U^n$.
\end{theorem}

Note that for lightlike hypersurfaces of the four-dimensional 
Minkowski space $M^4$ the result similar to the result of Theorem 
2 was obtained  in \cite{K89}. 

\section{The Fundamental Forms and Connections on a 
Lightlike Hypersurface of the de Sitter Space}

The first fundamental form of a lightlike hypersurface 
$U^n$ of the space $S^{n+1}_1$ is a metric quadratic form. It 
is defined by the scalar square of the differential $dx$ 
of a point of this  hypersurface. Since we have $x = A_n$, 
by (16) and (2) this scalar square has the form 
\begin{equation}\label{eq:22}
(dA_n, dA_n) = g_{ij} \omega_n^i  \omega_n^j = g
\end{equation}
and is a positive semidefinite 
differential quadratic form of signature 
$(n-1, 0)$. It follows that the system of equations 
$ \omega_n^i = 0$ defines on the  hypersurface  $U^n$ 
a fibration of isotropic lines which, as we showed in Section 2, 
coincide with rectilinear generators of this hypersurface.

The second fundamental form of a lightlike hypersurface 
$U^n$ determines its deviation from the tangent hyperplane 
$\eta$. To find this quadratic form, we compute the part of 
the second differential of the point $A_n$ which 
does not belong to the tangent hyperplane $\eta = A_0 \wedge A_1
\wedge \ldots \wedge A_n$:
$$
d^2 A_n \equiv \omega_n^i \omega_i^{n+1} A_{n+1} \pmod{\eta}.
$$
This implies that the second fundamental form can be written as 
\begin{equation}\label{eq:23}
b = \omega_n^i  \omega_i^{n+1} = \nu_{ij} \omega_n^i  \omega_n^j,
\end{equation}
where we used expression (20) for the form $\omega_i^{n+1}$. 
Since we assumed that $\mbox{{\rm rank}} \;
(\nu_{ij}) = n - 1$, the rank of the quadratic form (23) 
as well as the rank of the form $g$ is equal to $n - 1$. 
The nullspace of this   quadratic form  (see \cite{ON83}, p. 53) 
is again determined by  the system of equations 
$ \omega_n^i = 0$  and coincides with the isotropic direction 
on the hypersurface $U^n$. The reduction of the rank of the 
quadratic form $b$ is connected with the tangential degeneracy 
of the hypersurface $U^n$. The latter was noted in Theorem 2.

On a hypersurface $V^{n-1}$ of the conformal space $C^n$ 
that corresponds to a lightlike hypersurface $U^n \subset 
S^{n+1}_1$, the quadratic forms (22) and (23) define 
the net of curvature lines, that is, an orthogonal 
and conjugate net.

To find the connection forms of the hypersurface $U^n$, 
we find exterior derivatives of its basis forms 
$\omega_n^0$ and $\omega_n^i$:
\begin{equation}\label{eq:24}
\renewcommand{\arraystretch}{1.3}
\left\{
\begin{array}{ll}
d \omega_n^0 =  \omega_n^0 \wedge \omega_0^0 
                  + \omega_n^i \wedge \omega_i^0, \\
d \omega_n^i =  \omega_n^0 \wedge \omega_0^i 
                 + \omega_n^j \wedge \omega_j^i.
\end{array}
\renewcommand{\arraystretch}{1}
\right.
\end{equation}
This implies that the matrix 1-form 
\begin{equation}\label{eq:25}
\omega = \pmatrix{\omega_0^0 & \omega_i^0 \cr
                   \omega_0^i & \omega_j^i}
\end{equation} 
defines a torsion-free connection 
on the hypersurface $U^n$. To clarify the properties 
 of this connection, we find its curvature forms. To this end, we 
 substitute the null value $\omega_n^{n+1} = 0$ of 
the form $\omega_n^{n+1}$,  which by (17) 
defines $U^n$ along with the frame subbundle associated with 
$U^n$ in the space $S_1^{n+1}$, into expression (12) for  
the curvature forms of the de Sitter space $S_1^{n+1}$. As a 
result, we obtain 
\begin{equation}\label{eq:26}
\renewcommand{\arraystretch}{1.3}
\left\{
\begin{array}{lll}
\Omega_0^0 = d \omega_0^0 - \omega_0^i \wedge \omega_i^0 = 0, \\
\Omega_0^i = d \omega_0^i - \omega_0^0 \wedge \omega_0^i 
                   -  \omega_0^j \wedge \omega_j^i = 0, \\
\Omega_i^0 = d \omega_i^0 - \omega_i^0 \wedge \omega_0^0 
                   -  \omega_i^j \wedge \omega_j^0 
= -g_{ij}  \omega_n^j \wedge \omega_n^0, \\
\Omega_j^i = d \omega_j^i - \omega_j^0 \wedge \omega_0^i 
                   -  \omega_j^k \wedge \omega_k^i  
                   -  \omega_j^{n+1} \wedge \omega_{n+1}^i  
= -g_{jk}  \omega_n^k \wedge \omega_n^i. 
\end{array}
\renewcommand{\arraystretch}{1}
\right.
\end{equation}
In these formulas the forms $\omega_j^{n+1}$ and $\omega_0^i$ 
are expressed in terms of the basis forms $\omega_n^i$, and the 
forms $\omega_0^j, \omega_j^i$, and $\omega_i^0$ are 
fiber forms. If the principal parameters are fixed, then 
these fiber forms are invariant forms of the group $G$ 
of admissible transformations of frames associated with 
a point $x = A_n$ of the hypersurface $U^n$, and the connection 
 defined by the form (25) is a $G$-connection. 

To assign an affine connection on the hypersurface $U^n$, 
it is necessary to make a reduction of 
the family of frames in such a way that 
the forms $\omega_i^0$ become principal. Denote by 
$\delta$ the symbol of differentiation with respect to 
the fiber parameters, that is, for a fixed point $x = A_n$ 
of the hypersurface $U^n$, and by $\pi_\eta^\xi$ the values 
of the 1-forms $\omega_\eta^\xi$ for a fixed point $x = A_n$, 
that is, $\pi_\eta^\xi = \omega_\eta^\xi (\delta)$. Then 
we obtain
$$
\pi_n^0 = 0, \pi_n^i = 0, \pi_i^n = 0, \pi_i^{n+1} = 0.
$$
It follows 
\begin{equation}\label{eq:27}
\delta A_i = \pi_i^0 A_0 + \pi_i^j A_j.
\end{equation}

The points $A_0$ and $A_i$ determine the tangent subspace to the 
submanifold $V^{n-1}$ described by the point $A_0$ on the
 hyperquadric $Q^n$. If we fix an $(n-2)$-dimensional subspace 
$\zeta$ not containing the point $A_0$ in this tangent subspace 
and place the points $A_i$ into $\zeta$, then we obtain $\pi_i^0$. 
This means that the forms 
$\omega_i^0$ become principal, that is, 
\begin{equation}\label{eq:28}
\omega_i^0 = \mu_{ij} \omega_n^j + \mu_i \omega_n^0,
\end{equation}
and as a result, an affine connection arises 
on the hypersurface $U^n$. 

We will call the subspace $\zeta \subset T_{A_0} (V^{n-1})$ 
the {\em normalizing subspace} of the lightlike 
hypersurface $U^n$. We have proved the following result:

\begin{theorem} If in every tangent subspace $T_{A_0} (V^{n-1})$ 
of  the submanifold $V^{n-1}$ associated with a lightlike 
hypersurface $U^n, V^{n-1} \subset Q^n$, a normalizing 
 $(n-2)$-dimensional subspace $\zeta$  not containing the point 
$A_0$ is assigned, then there arises a torsion-free 
affine connection on $U^n$. 
\end{theorem}

The last statement of Theorem 4 follows from the first 
two equations of (26). 

By (28), the last equation of (26) can be written in the form
\begin{equation}\label{eq:29}
\renewcommand{\arraystretch}{1.3}
\begin{array}{ll}
\widetilde{\Omega}_j^i &= d \omega_j^i 
                   -  \omega_j^k \wedge \omega_k^i  
= g^{im}(-g_{jk} g_{ml} + \mu_{jk} \nu_{ml} + 
\nu_{jk} \mu_{ml}) \omega_n^k \wedge \omega_n^l \\
&+ g^{im}(\mu_j \nu_{ml} - \mu_m \nu_{jl}) 
\omega_n^0 \wedge \omega_n^l. 
\end{array}
\renewcommand{\arraystretch}{1}
\end{equation}
From the first three equations of (26) and equation (29) 
we can find the torsion tensor of the affine connection 
indicated in Theorem 4:
\begin{equation}\label{eq:30}
\renewcommand{\arraystretch}{1.3}
\begin{array}{ll}
R^0_{0uv} = 0, \;\; R^i_{0uv} = 0, 
\;\; R^0_{ij0} = - R^0_{i0j} = - \frac{1}{2} g^{ij}, \\
R^i_{jkl} =  \frac{1}{2} g^{im} (g_{jl} g_{mk} - g_{jk} g_{ml} 
+ \mu_{jk} \nu_{ml} -  \mu_{jl} \nu_{mk} \\
\;\;\;\;\;\;\;\;\;\;\; 
+  \nu_{jk} \mu_{ml} -  \nu_{jl} \mu_{mk}), \\
R^i_{j0l} = - R^i_{jl0} = \frac{1}{2} g^{im} (\mu_{j} \nu_{ml} - \mu_{m} \nu_{jl}).
\end{array}
\renewcommand{\arraystretch}{1}
\end{equation}

The constructed above fibration of normalizing subspaces $\zeta$ 
defines a distribution $\Delta$ of $(n-1)$-dimensional 
elements on a lightlike hypersurface $U^n$. In fact, 
the point $x = A_n$ of the hypersurface $U^n$ along with the 
subspace $\zeta =  A_1 \wedge \ldots \wedge 
A_{n-1}$ define the $(n-1)$-dimensional subspace 
which is complementary to the straight line $A_n A_0$ and 
lies in the tangent subspace $\eta$ of the hypersurface $U^n$. 
Following the book \cite{DB96}, we will call this subspace 
the {\em screen}, and the distribution $\Delta$ the 
{\em screen distribution}. Since at the point $x$ the screen 
is determined by the subspace $A_n A_1 \ldots A_{n-1}$, 
the differential equations of the screen distribution 
has the form
\begin{equation}\label{eq:31}
\omega_n^0 = 0.
\end{equation}
But by (28)
$$
d\omega_n^0 = \omega_n^i \wedge (\mu_{ij} \omega_n^j 
+ \mu_i \omega_n^0).
$$
Hence the screen distribution  is integrable if and only if the tensor $\mu_{ij}$  is symmetric. Thus we arrived 
at the following result:

\begin{theorem} The fibration of normalizing subspaces $\zeta$ 
defines a screen distribution $\Delta$ of $(n-1)$-dimensional 
elements on a lightlike hypersurface $U^n$. This  distribution 
is integrable if and only if the tensor $\mu_{ij}$ 
defined by equation $(28)$ is symmetric.
\end{theorem}

Note that the configurations similar to that described in Theorem 
5 occurred in the works of the Moscow geometers published in the 
1950s. They were called the {\em one-side stratifiable pairs of 
ruled surfaces} (see \cite{F56}, \$30 or \cite{AG93}, p. 187).

\section{An Invariant Normalization of 
Lightlike 
\\
Hypersurfaces of the de Sitter Space}

In \cite{A52} (see also \cite{AG96}, Ch. 2) 
an invariant normalization of a hypersurfaces $V^{n-1}$ 
of the conformal space $C^n$ was constructed. By Theorem 3, this normalization can be interpreted in terms of the geometry of 
the de Sitter space $S_1^{n+1}$.

Taking exterior derivative of equations (18) defining 
the hypersurface $V^{n-1}$ in the conformal space $C^n$, 
we obtain
$$
 \omega_i^n \wedge  \omega_0^i = 0,
$$
from which by linear independence of the 1-forms 
$\omega_0^i$ on $V^{n-1}$ and Cartan's lemma we find that 
\begin{equation}\label{eq:32}
\omega_i^n = \lambda_{ij} \omega_0^j, \;\; \lambda_{ij} 
= \lambda_{ji}.
\end{equation}
Here and in what follows we retain the notations used 
in the study of the geometry of hypersurfaces  
of the conformal space $C^n$ in the book \cite{AG96}. 

It is not difficult to find relations between the 
coefficients $\nu_{ij}$ in formulas (20) and 
 $\lambda_{ij}$ in formulas (32). Substituting the values 
of the forms $\omega_i^n$ and $\omega_0^j$ from (5) into 
(32), we find that 
$$
- g_{ij} \omega_n^j = \lambda_{ij} g^{jk} \omega_k^{n+1}.
$$
Solving these equations for $\omega_k^{n+1}$, we obtain 
$$
 \omega_i^{n+1} = - g_{ik} \widetilde{\lambda}^{kl} g_{lj}  \omega_n^j, 
$$
where $(\widetilde{\lambda}^{kl})$ is the inverse matrix of the 
matrix $(\lambda_{ij})$. Comparing these equations with equations (20), we obtain
\begin{equation}\label{eq:33}
\nu_{ij} = - g_{ik} \widetilde{\lambda}^{kl} g_{lj}.
\end{equation}
Of course, in this computation we assumed that the 
matrix $(\lambda_{ij})$ is nondegenerate.

Let us clarify the geometric meaning of the vanishing of 
$\det (\lambda_{ij})$. To this end, we make an admissible 
transformation of the moving frame associated with a point 
of a lightlike hypersurface $U^n$ by setting 
\begin{equation}\label{eq:34}
\widehat{A}_n = A_n + s A_0.
\end{equation}
The point $\widehat{A}_n$ as the point $A_n$ lies on the 
rectilinear generator $A_n A_0$. Differentiating this point and 
applying formulas (16) and (19), we obtain
\begin{equation}\label{eq:35}
d\widehat{A}_n = (ds + s\omega_0^0 + \omega_n^0) A_0 + 
(\omega_n^i + s \omega_0^i) A_i.
\end{equation}
It follows that in the new frame the form $\omega_n^i$ 
becomes 
$$
\widehat{\omega}_n^i = \omega_n^i + s \omega_0^i.
$$
By (5) and (32), it follows that 
$$
\widehat{\omega}_n^i = -g^{ik} (\lambda_{kj} - s g_{kj}) \omega_0^j.
$$
This implies that in the new frame the quantities $\lambda_{ij}$ 
become
\begin{equation}\label{eq:36}
\widehat{\lambda}_{ij} = \lambda_{ij} - s g_{ij}.
\end{equation}
Consider also the matrix $(\widehat{\lambda}^i_j) 
= (g^{ik} \widehat{\lambda}_{kj})$. Since $g_{ij}$ is 
a nondegenerate tensor, the matrices  $(\widehat{\lambda}^i_j)$ 
and $(\widehat{\lambda}_{ij})$ have the same rank 
$\rho \leq n - 1$. 

From equation (35) it follows that 
$$
d\widehat{A}_n = (ds + s\omega_0^0 + \omega_n^0) A_0  
- \widehat{\lambda}^i_j A_i \omega^j_0.
$$
Hence the tangent subspace to 
the  hypersurface $U^n$ at the point $\widehat{A}_n$ 
is determined by the points $\widehat{A}_n, A_0$, and 
$\widehat{\lambda}^i_j A_i$. At the points,  at which 
the rank $\rho$ of the 
matrix $(\widehat{\lambda}^i_j)$ is equal to $n-1, \rho 
= n - 1$, the 
tangent subspace to the  hypersurface $U^n$ has dimension $n$, 
and such points are {\em regular points} of the  hypersurface. 
The points,  at which the rank $\rho $ of the 
matrix $(\widehat{\lambda}^i_j)$ is reduced, 
are {\em singular points} of the  hypersurface $U^n$. The 
coordinates 
of singular points are defined by the condition 
$\det (\widehat{\lambda}^i_j) = 0$ which by (36) is equivalent to 
the equation 
\begin{equation}\label{eq:37}
\det (\lambda_{ij} - s g_{ij}) = 0,
\end{equation}
the {\em characteristic equation} of the matrix $(\lambda_{ij})$ with respect to the tensor $g_{ij}$. The degree of this equation 
is equal to $n - 1$.

In particular, if  $A_n$ is a regular point 
of the  hypersurface $U^n$, then the matrix $(\lambda_{ij})$ 
is nondegenerate, and equation (33) holds. On the other hand, 
if  $A_n$ is a singular point 
of  $U^n$, then  equation (33) is meaningless. 

Since the matrix $(\lambda_{ij})$ is symmetric and 
the matrix $(g_{ij})$ defines a positive definite form 
of rank $n - 1$, equation (37) has 
$n - 1$ real roots if each root is counted as many times 
as its multiplicity. Thus on a rectilinear generator $A_n A_0$ of 
a lightlike   hypersurface $U^n$ there are $n-1$ real singular points. 

By Vieta's theorem, the sum of the roots of equation (37) is 
equal to the coefficient in $s^{n-2}$, and this coefficient 
is $\lambda_{ij} g^{ij}$. Consider the quantity
\begin{equation}\label{eq:38}
\lambda = \frac{1}{n-1} \lambda_{ij} g^{ij},
\end{equation}
which is the arithmetic mean of the roots of equation (37). 
This quantity $\lambda$ allows us to construct  new 
quantities
\begin{equation}\label{eq:39}
a_{ij} = \lambda_{ij} - \lambda g_{ij}.
\end{equation}
It is easy to check that the quantities $a_{ij}$ 
do not depend on the location of the point $A_n$ on the 
straight line $A_n A_0$, that is, $a_{ij}$ is invariant with 
respect to the transformation of the moving frame defined by 
equation (34). Thus the quantities $a_{ij}$ form a tensor 
on the  hypersurface $U^n$ defined in its 
second-order neighborhood. This tensor satisfies the condition 
\begin{equation}\label{eq:40}
a_{ij} g^{ij} = 0,
\end{equation}
that is, it is apolar to the tensor $g_{ij}$. 

On the straight line $A_n A_0$ we consider a point
\begin{equation}\label{eq:41}
C = A_n + \lambda A_0.
\end{equation}
It is not difficult to check that this point remains also 
fixed when the point $A_n$ moves along the  straight line $A_n A_0$. Since $\lambda$ is the arithmetic mean of the roots of equation (37) defining singular points on the  straight line $A_n A_0$, the point $C$ is the {\em harmonic pole} 
(see \cite{C50}) of the point $A_0$ with respect to these 
singular points. In particular, for $n = 3$, the point $C$ is the 
fourth harmonic point to the point $A_0$ with respect to two 
singular points of the rectilinear generator $A_3 A_0$ of the lightlike  hypersurface $U^3$ of the de Sitter space $S_1^4$. 

In the conformal theory of hypersurfaces, to 
the point $C$ there corresponds a hypersphere which is tangent 
to the hypersurface at the point $A_0$. This hypersphere
 is called the {\em central tangent hypersphere} 
(see \cite{AG96}, pp. 40--41). Since 
\begin{equation}\label{eq:42}
(d^2 A_0, C) = a_{ij} \omega_0^i \omega_0^j, 
\end{equation}
the cone 
$$
a_{ij} \omega_0^i \omega_0^j = 0
$$
with vertex at the point $A_0$ belonging to the tangent 
subspace $T_{A_0} (V^{n-1})$ contains the directions along which 
the central  hypersphere has a second-order tangency with 
the hypersurface $V^{n-1}$ at the point $A_0$. From the apolarity 
condition (39) it follows that it is possible to inscribe an orthogonal $(n-1)$-hedron with vertex at $A_0$ into the cone defined by equation (42). 

Now we can construct an invariant normalization of 
a lightlike hypersurface $U^n$ of the de Sitter space 
$S_1^{n+1}$. 
To this end, first we repeat some computations from Ch. 2 
of \cite{AG96}. 

Taking exterior derivatives of equations (32) and applying Cartan's lemma, we  obtain the equations
\begin{equation}\label{eq:43}
\nabla \lambda_{ij} +  \lambda_{ij} \omega_0^0 + g_{ij} \omega_n^0 =  \lambda_{ijk} \omega_0^k, 
\end{equation}
where 
$$
\nabla \lambda_{ij} = d \lambda_{ij} - \lambda_{ik} \omega^k_j 
-  \lambda_{kj} \omega^k_i, 
$$
and the quantities $ \lambda_{ijk}$ are symmetric 
with respect to all three indices. 
Equations (43) confirm one more time that the 
quantities $ \lambda_{ij} $ do not form a tensor 
and depend on a location of the point $A_n$ on the straight line 
$A_n A_0$. This dependence is described by a closed form 
relation (36). From formulas (43) it follows that the quantity 
$\lambda$ defined by equations (38) satisfy the differential 
equation
\begin{equation}\label{eq:44}
d \lambda  +  \lambda \omega_0^0 +  \omega_n^0 =  \lambda_k \omega_0^k, 
\end{equation}
where 
$$
\lambda_k = \frac{1}{n-1} g^{ij} \lambda_{ijk}
$$
(see formulas (2.1.35) and (2.1.36) in the book 
\cite{AG96}). 

The point $C$ lying on the rectilinear generator $A_n A_0$ 
of the  hypersurface $U^n$ describes a submanifold $W \subset U^n$ when  $A_n A_0$ moves. Let us find the tangent 
subspace to $U^n$ at the point $C$. Differentiating equation (40) 
and applying formulas (16) and (19), we obtain
$$
dC = (d \lambda + \lambda \omega_0^0 +  \omega_n^0) A_0 
+ (\omega_n^i +  \lambda \omega_0^i) A_i.
$$
By (5), (32), (39), and (44), it follows that 
 \begin{equation}\label{eq:45}
d C = (\lambda_i A_0  - g^{jk} a_{ki} A_j) \omega_0^i.   \end{equation}
Define the affinor
 \begin{equation}\label{eq:46}
a^i_j = g^{ik} a_{kj}, 
 \end{equation}
whose rank coincides with the rank of the tensor $a_{ij}$. 
Then equation (45) takes the form
$$
d C = (\lambda_i A_0  -  a_i^j A_j) \omega_0^i.  
$$
The points
 \begin{equation}\label{eq:47}
C_i = \lambda_i A_0  -  a_i^j A_j
 \end{equation}
together with the point $C$ define the tangent subspace 
to the submanifold $W$ described by the point $C$ on 
 the  hypersurface $U^n$. 

If the point $C$ is a regular point of 
the rectilinear generator $A_n A_0$ of the  hypersurface $U^n$, 
then  the rank of the tensor $a_{ij}$ defined by 
equations (39) as well as the rank of the affinor $a_j^i$ is 
equal to $n-1$. As a result, the points $C_i$ are linearly 
independent and together with the point $C$ define 
the $(n-1)$-dimensional tangent subspace $T_C (W)$, and the submanifold $W$ itself has dimension $n-1, \dim W = n-1$.

The points $C_i$ also belong to the tangent subspace 
$T_{A_0} (V^{n-1})$ and define  the $(n-2)$-dimensional 
subspace $\zeta = T_{A_0} (V^{n-1}) \cap T_C (W)$ in it. This 
subspace 
is a normalizing subspace. Since such a normalizing subspace is 
defined in each tangent subspace $T_{A_0} (V^{n-1})$ of the 
hypersurface $V^{n-1} \subset Q^n$, there arises the fibration of 
these subspaces which by Theorem 4 defines an invariant affine 
connection on the lightlike hypersurface $U^n$.
 Thus we proved the following result:

\begin{theorem} If the tensor $a_{ij}$ defined by formula $(39)$ 
 on a lightlike hypersurface $U^n \subset S^{n+1}_1$ 
is nondegenerate, 
then it is possible to construct the invariant normalization 
of $U^n$ by means of the $(n-2)$-dimensional subspaces
$$
\zeta = C_1 \wedge C_2 \wedge \ldots \wedge C_{n-1}.
$$
This normalization induces on $U^n$ an invariant screen 
distribution and an invariant affine 
connection intrinsically connected with the geometry of this 
hypersurface.
\end{theorem}

Theorem 5 implies that the invariant normalization 
we have constructed defines on $U^n$ an invariant screen 
distribution $\Delta$ which is also 
intrinsically connected with the geometry of the 
hypersurface $U^n$; here $\Delta_x = x \wedge \xi, x \in A_n A_0$.

Note that for the hypersurface $V^{n-1}$ of the conformal space $C^n$ 
a similar  invariant normalization was constructed as far 
back as 1952 (see \cite{A52} and also \cite{AG96}, Ch. 2). 
In the present paper we gave a new geometric meaning of 
this invariant normalization. 

\section{Singular Points of Lightlike  Hypersurfaces of 
the de Sitter Space}

As we indicated in Section 4, the points 
\begin{equation}\label{eq:48}
z =  A_n  + s A_0
\end{equation}
of the rectilinear generator $A_n A_0$ of the lightlike 
hypersurface $U^n$ are singular if their nonhomogeneous 
coordinate $s$ satisfies the equation 
\begin{equation}\label{eq:49}
\det (\lambda_{ij} - s g_{ij}) = 0.
\end{equation}
We will investigate in more detail the structure of 
a lightlike hypersurface $U^n$ in a neighborhood of its singular 
point. 

Equation (49) is the characteristic equation of 
the matrix $(\lambda_{ij})$ with respect to the tensor 
$(g_{ij})$. The degree of this equation is $n-1$, and 
since the matrix  $(\lambda_{ij})$  is symmetric and 
 the matrix $(g_{ij})$ is also symmetric and 
positive definite, then according 
to the well-known result of linear algebra, all  roots 
of this equation are real, and the matrices $(\lambda_{ij})$  
 and  $(g_{ij})$ can be simultaneously reduced 
to a diagonal form. 

 Denote the roots of the characteristic equation 
 by $s_h, h = 1, 2, 
\ldots , n - 1$, and denote the corresponding singular points 
of the rectilinear generator $A_n A_0$ by 
\begin{equation}\label{eq:50}
B_h =  A_n  + s_h A_0.
\end{equation}
These singular points are called {\em foci}
 of the rectilinear generator $A_n A_0$ of a lightlike 
hypersurface $U^n$. 

It is clear from (50) that the point $A_0$ is not a focus 
 of the rectilinear generator $A_n A_0$. This is explained 
by the fact that by our assumption $\mbox{{\rm rank}} \; 
(\nu_{ij}) = n - 1$, and by (21), on the hyperquadric $Q^n$ 
the point $A_0$ describes a hypersurface $V^{n-1}$ which is 
transversal to the straight lines $A_0 A_n$. 

In the conformal theory of hypersurfaces, to 
the singular points $B_h$ there correspond the 
tangent hyperspheres defining the principal directions at 
a point $A_0$ of the hypersurface $V^{n-1}$ of the conformal 
space $C^n$ (see \cite{AG96}, p. 55). 

We will construct a classification of singular points 
of a lightlike hypersurface $U^n$ of the space $S^{n+1}_1$. 
We will use some computations that we made while constructing 
a classification of canal hypersurfaces in \cite{AG97}. 

Suppose first that $B_1 = A_n + s_1 A_0$ be a singular 
point defined by a simple 
root $s_1$ of characteristic equation (49), 
$s_1 \neq s_h, h = 2, \ldots , n - 1$. For this singular point 
we have 
\begin{equation}\label{eq:51}
dB_1 = (ds_1 + s_1 \omega_0^0 + \omega_n^0) A_0 - 
\widehat{\lambda}_j^i \omega_0^j A_i,
\end{equation}
where 
\begin{equation}\label{eq:52}
\widehat{\lambda}_j^i = g^{ik} (\lambda_{kj} - s_1 g_{kj}) 
\end{equation}
is a degenerate symmetric affinor having a single null 
eigenvalue. The matrix of this affinor can be reduced 
to a quasidiagonal form
\begin{equation}\label{eq:53}
(\widehat{\lambda}_j^i) = \pmatrix{0 & 0 \cr 
             0 & \widehat{\lambda}_q^p}, 
\end{equation}
where $p, q = 2, \ldots , n - 1$, and $(\widehat{\lambda}_q^p)$ 
is a nondegenerate symmetric affinor. The matrices 
$(g_{ij})$ and $(\lambda_{ij} - s_1 g_{ij})$ are reduced 
to the forms
$$
\pmatrix{1 & 0 \cr 
         0 & g_{pq}} \;\; \mbox{{\rm and}} \;\;
\pmatrix{0 & 0 \cr 
             0 & \widehat{\lambda}_{pq}}, 
$$
where $(\widehat{\lambda}_{pq}) = (\lambda_{pq} - s_1 
g_{pq})$ is a nondegenerate symmetric matrix. 

Since the point $B_1$ is defined invariantly on the 
generator $A_n A_0$, then it will be fixed if $\omega_0^i = 0$. 
Thus it follows from (51) that
\begin{equation}\label{eq:54}
ds_1 + s_1 \omega_0^0 + \omega_n^0 = s_{1i} \omega^i,
\end{equation}
here and in what follows $\omega^i = \omega^i_0$. 
By (53) and (54) relation (51) takes the form 
\begin{equation}\label{eq:55}
dB_1 = s_{11} \omega^1 A_0 + (s_{1p} A_0 
- \widehat{\lambda}_p^q A_q) 
 \omega^p.
\end{equation}
Here the points $C_p = s_{1p} A_0 - \lambda_p^q A_q$ 
are linearly independent and belong to the tangent subspace 
$T_{A_0} (V^{n-1})$. 

Consider the submanifold ${\cal F}_1$ 
described by the singular point 
$B_1$ in the space $S^{n+1}_1$. This submanifold is called 
the {\em focal manifold} of the hypersurface $U^n$. 
Relation (55) shows that two cases are possible:

\begin{description}
\item[1)] $s_{11} \neq 0$. In this case  the submanifold 
${\cal F}_1$ 
is of dimension $n - 1$, and its tangent subspace at the point 
$B_1$ is determined by the points $B_1, A_0$, and $C_p$. This 
subspace contains the straight line $A_n A_0$, intersects 
the hyperquadric $Q^n$, and thus it,  as well as  the submanifold 
${\cal F}_1$ itself, is timelike. For $\omega^p = 0$, 
the point $B_1$ describes a curve $\gamma$ on the submanifold 
${\cal F}_1$  
which is tangent to the straight line $B_1 A_0$ coinciding with 
the generator $A_n A_0$ of the hypersurface $U^n$. 
The curve $\gamma$ is an isotropic curve of the de Sitter 
space $S^{n+1}_1$. Thus on 
${\cal F}_1$ there arises a fibre bundle of focal lines. The 
hypersurface $U^n$ is foliated into an $(n-2)$-parameter 
family of torses for which these lines are edges of regressions. 
The points $B_1$ are singular points of a kind which is called 
a {\em fold}. 

If the characteristic equation (49) has distinct roots, then 
an isotropic rectilinear generator $l$ of a lightlike 
hypersurface 
$U^n$ carries $n - 1$ distinct foci $B_h, h = 1, \ldots , n-1$. 
If for each of these foci the condition of type $s_{11} \neq 0$ holds, 
then each of them  describes a focal submanifold ${\cal F}_h$,  
carrying conjugate net.  Curves of one family of this net are tangent 
to the straight lines $l$, and this family is isotropic. 
On the hypersurface $V^{n-1}$ of the space $C^n = Q^n$ described by 
the point $A_0$, to these conjugate nets there correspond the 
net of curvature lines. 

\item[2)] $s_{11} = 0$. In this case  relation (55) takes 
the form
\begin{equation}\label{eq:56}
dB_1 =  (s_{1p} A_0 - \widehat{\lambda}_p^q A_q)  \omega^p,
\end{equation}
and the focal submanifold ${\cal F}_1$ 
is of dimension $n - 2$. Its tangent subspace at the point 
$B_1$ is determined by the points $B_1$ and $C_p$. 
An arbitrary point $z$ of this subspace can be written in the 
form 
$$
z = z^n B_1 + z^p C_p = z^n (A_n + s_1 A_0) + z^p 
(s_{1p} A_0 - \widehat{\lambda}_p^qA_q).
$$
Substituting the coordinates of this point into relation (3), 
we find that 
$$
(z, z) = g_{rs} \widehat{\lambda}_p^r   
\widehat{\lambda}_q^s z^p z^q + (z^n)^2 > 0.
$$
It follows that the tangent subspace 
$T_{B_1} (F_1)$ does not have common points with 
the hyperquadric $Q^n$, that is, it is spacelike. Since this 
takes place for any point $B_1 \in {\cal F}_1$, the focal 
submanifold ${\cal F}_1$ is spacelike.

For $\omega^p = 0$, the point $B_1$ is fixed. The subspace 
$T_{B_1} ({\cal F}_1)$ will be fixed too. 
On the hyperquadric $Q^n$, the point $A_0$ describes 
a curve $q$ which is polar-conjugate to $T_{B_1} ({\cal F}_1)$. 
Since $\dim T_{B_1} ({\cal F}_1) = n - 2$, the curve $q$ is 
a conic, along which the two-dimensional plane 
polar-conjugate to the subspace $T_{B_1} ({\cal F}_1)$  with 
respect to the hyperquadric $Q^n$, intersects $Q^n$. 
Thus for $\omega^p = 0$, the rectilinear generator $A_n A_0$ of 
the hypersurface $U^n$ describes  a two-dimensional 
second-order cone  with vertex at the point $B_1$ 
and the directrix $q$. Hence in the case under consideration 
a lightlike hypersurface $U^n$ is foliated into an 
$(n-2)$-parameter family of second-order cones whose vertices 
describe the $(n-2)$-dimensional focal submanifold 
${\cal F}_1$, and 
the points $B_1$ are {\em conic} singular points of 
the hypersurface $U^n$. 

The hypersurface $V^{n-1}$ of the conformal space $C^n$ 
corresponding to such a lightlike hypersurface $U^n$ 
is a canal hypersurface which envelops an $(n-2)$-parameter 
family of hyperspheres. Such a hypersurface carries a family 
of cyclic generators which depends on the same number 
of parameters. Such hypersurfaces were investigated in detail 
in \cite{AG97}. 
\end{description}

Further let $B_1$ be a  singular point of multiplicity $m$, where 
$m \geq 2$, of a rectilinear generator $A_n A_0$ 
of  a lightlike hypersurface $U^n$ of the space $S^{n+1}_1$ 
defined by an $m$-multiple root of characteristic equation 
(49). We will assume that 
\begin{equation}\label{eq:57}
 s_1 = s_2 = \ldots = s_m := s_0, s_0 \neq s_p, 
\end{equation}
and also assume that $a, b, c = 1, \ldots, m$ 
and $p, q, r = m+1, \ldots, n-1$. Then the matrices $(g_{ij})$ 
and $(\lambda_{ij})$ can be simultaneously  
reduced to  quasidiagonal forms
$$
\pmatrix{g_{ab} & 0 \cr 
             0 & g_{pq}} \;\; \mbox{{\rm and}} \;\;
\pmatrix{s_0 g_{ab} & 0 \cr 
             0 & \lambda_{pq}}.
$$
We also construct the matrix $(\widehat{\lambda}_{ij}) 
= (\lambda_{ij} - s_0 g_{ij})$. Then 
\begin{equation}\label{eq:58}
(\widehat{\lambda}_{ij}) = \pmatrix{0 & 0 \cr 
             0 & \widehat{\lambda}_{pq}},
\end{equation}
where $\widehat{\lambda}_{pq} = \lambda_{pq} - s_0 g_{pq}$ 
is a nondegenerate matrix of order $n - m - 1$. 

By relations (58) and formulas (5) and (32) we have
\begin{equation}\label{eq:59}
\omega_a^n - s_0 \omega_a^{n+1} = 0,
\end{equation}
\begin{equation}\label{eq:60}
\omega_p^n - s_0 \omega_p^{n+1} =  \widehat{\lambda}_{pq} 
\omega^q.
\end{equation}
Taking exterior derivative of equation (59) and applying relation 
(60), we find that 
\begin{equation}\label{eq:61}
\widehat{\lambda}_{pq} \omega_a^p \wedge \omega^q 
+ g_{ab} \omega^b \wedge (ds_0 + s_0 \omega_0^0 
+ \omega_n^0) = 0.
\end{equation}
It follows that the 1-form $ds_0 + s_0 \omega_0^0 
+ \omega_n^0$ can be expressed in terms of the basis forms. 
We write these expressions in the form
\begin{equation}\label{eq:62}
ds_0 + s_0 \omega_0^0 + \omega_n^0 = s_{0c} \omega^c + 
s_{0q} \omega^q.
\end{equation}
Substituting this decomposition into equation (61), we find that 
\begin{equation}\label{eq:63}
(\widehat{\lambda}_{pq} \omega_a^p  
+ g_{ab} s_{0q} \omega^b) \wedge \omega^q + g_{ab}  s_{0c} 
\omega^b \wedge \omega^c  = 0.
\end{equation}
The terms in the left hand side of (63) do not have similar 
terms. Hence both terms are equal to 0. Equating to 0 
the coefficients of the summands of the second term, 
we find that 
\begin{equation}\label{eq:64}
 g_{ab}  s_{0c} =  g_{ac}  s_{0b}.
\end{equation}
Contracting this equation with the matrix $(g^{ab})$ 
which is the inverse matrix of the matrix $(g_{ab})$, we obtain
$$
ms_{0c} = s_{0c}.
$$
Since $m \geq 2$, it follows that 
$$
s_{0c} = 0,
$$
and  relation (62) takes the form 
\begin{equation}\label{eq:65}
ds_0 + s_0 \omega_0^0 + \omega_n^0 = s_{0p} \omega^p.
\end{equation}

For the  singular point of multiplicity $m$ of the generator 
$A_n A_0$ in question the equation (51) can be written in the 
form
$$
dB_1 = (ds_0 + s_0 \omega_0^0 + \omega_n^0) A_0 - 
\widehat{\lambda}_q^p \omega_0^q A_p.
$$
Substituting decomposition (65) in the last equation, 
we find that 
\begin{equation}\label{eq:66}
dB_1 = (s_{0p} A_0 - \widehat{\lambda}_p^q A_q)\omega_0^p.
\end{equation}
This relation is similar to equation (56) with the only 
difference that in (56) we had $p, q = 2, \ldots, n - 1$, and 
in (66) we have $p, q = m+1, \ldots, n - 1$. 
Thus the point $B_1$ describes now a spacelike focal 
manifold ${\cal F}_1$ of dimension $n-m-1$. For $\omega_0^p = 0$, 
the point $B_1$ is fixed, and the point $A_0$ describes an 
$m$-dimensional submanifold on the hyperquadric $Q^n$ which 
is a cross-section of $Q^n$ by an $(m+1)$-dimensional subspace 
that is polar-conjugate to the $(n-m-1)$-dimensional subspace 
tangent to the submanifold ${\cal F}_1$. 

The point $B_1$ is a  conic singular point of multiplicity $m$ of 
a lightlike hypersurface $U^n$, and this hypersurface is foliated 
into an $(n-m-1)$-parameter family of $(m+1)$-dimensional 
second-order cones circumscribed about 
the hyperquadric $Q^n$. The  hypersurface $V^{n-1}$ of the 
conformal space $C^n$ that corresponds to such a  hypersurface 
$U^n$ is an $m$-canal hypersurface (i.e.,  
the envelope of an $(n-m-1)$-parameter family of 
hyperspheres), and it carries an $m$-dimensional 
spherical generators.

Note also an extreme case when the rectilinear generator 
$A_n A_0$ of a lightlike hypersurface $U^n$ carries a single 
 singular point of multiplicity $n-1$. 
As follows from our consideration 
of the cases $m \geq 2$, this singular point is fixed, and 
the hypersurface $U^n$ become a second-order hypercone with 
vertex at this singular point which is circumscribed about 
 the hyperquadric $Q^n$. This hypercone is the isotropic cone 
of the space $S^{n+1}_1$. The  hypersurface $V^{n-1}$ of the 
conformal space $C^n$ that corresponds to such a  hypersurface 
$U^n$ is a hypersphere of the space $C^n$.
 
The following theorem combines the results of this section:

\begin{theorem}
 A lightlike hypersurface $U^n$ of maximal rank $r=n-1$ of the de 
Sitter space $S^{n+1}_1$ possesses $n-1$ real singular points 
on each of its rectilinear generators if each of these singular 
points  is counted as many times as its multiplicity. The simple 
singular points can be of two kinds: a fold and conic. In the 
first case the hypersurface $U^n$ is foliated into an $(n-2)$-parameter family of torses, and in the second case it is foliated 
into an $(n-2)$-parameter family of second-order cones. 
The vertices of these cones 
describe the $(n-2)$-dimensional spacelike submanifold 
in the space $S^{n+1}_1$. All multiple singular points 
of a hypersurface $U^n$ are conic. If a rectilinear 
generator of a hypersurface $U^n$ carries a singular point 
of multiplicity $m$, $2 \leq m \leq n -1$, 
then the hypersurface $U^n$ is foliated into an $(n-m-1)$-parameter family of $(m+1)$-dimensional second-order cones. 
The vertices of these cones 
describe the $(n-m-1)$-dimensional spacelike submanifold 
in the space $S^{n+1}_1$.
The hypersurface $V^{n-1}$ of the conformal space $C^n$ 
corresponding to a lightlike hypersurface $U^n$ with 
 singular points of multiplicity $m$ 
is a canal hypersurface which envelops an $(n-m-1)$-parameter 
family of hyperspheres and has $m$-dimensional spherical 
generators. 
\end{theorem}

Since lightlike hypersurfaces $U^n$ of 
the de Sitter space $S^{n+1}_1$ represent a light flux 
(see Section {\bf 2}), its focal submanifolds have the following 
physical meaning. If one of them is a lighting submanifold, then 
others will be manifolds of concentration of a light flux. 
Intensity of concentration depends on multiplicity of a focus 
 describing this submanifold. 

In the extreme case when an isotropic rectilinear generator 
$l = A_n A_0$ of a hypersurface $U^n$ carries one $(n-1)$-multiple 
focus, the hypersurfaces $U^n$ degenerates into the light cone 
generated by a point source of light. This cone represents a radiating 
light flux.

If each isotropic generator $l \subset U^n$ carries two foci 
$B_1$ and $B_2$ of multiplicities $m_1$ and $m_2, \,m_1 + m_2 = n - 1, \,
m_1 > 1, \, m_2 > 1$, then these foci describe spacelike submanifolds 
${\cal F}_1$ and ${\cal F}_2$ of dimension $n - m_1 - 1$ and 
$n - m_2 - 1$, respectively. If one of these submanifolds is a lighting 
 submanifold, then on the second one a light flux is concentrated.

\section{Lightlike  Hypersurfaces of Reduced Rank}

As we proved in Section 2, lightlike hypersurfaces 
of the de Sitter space $S^{n+1}_1$ are ruled tangentially 
degenerate   hypersurfaces. However in all preceding sections 
starting from Section 3 we assumed that the rank of these 
  hypersurfaces is maximal, that is, it is equal to $n - 1$. 
In this section we consider lightlike  hypersurfaces 
of reduced rank $r < n - 1$.

We proved in Section 2 that the rank of a lightlike  hypersurface 
$U^n$ coincides with the rank of the matrix $(\nu_{ij})$ 
defined by equation (20) as well as with the dimension of the 
submanifold $V$ described by the point $A_0$ on the Darboux 
hyperquadric $Q^n$. As a result, to a lightlike  hypersurface 
$U^n$ of rank $r$ there corresponds an $r$-dimensional 
submanifold $V = V^r$ in the conformal space $C^n$.

The symmetric matrices $(g_{ij})$ and $(\nu_{ij})$ first of which 
is nondegenerate and positive definite and second is of rank 
$r$, can be simultaneously reduced to  quasidiagonal forms
\begin{equation}\label{eq:67}
(g_{ij}) = \pmatrix{g_{ab} & 0 \cr 
             0 & g_{pq}} \;\; \mbox{{\rm and}} \;\;
(\nu_{ij}) = \pmatrix{0  & 0 \cr 
             0 & \nu_{pq}},
\end{equation}
where $a, b = 1, \ldots , m;\, p, q, s = m+1, \ldots , n - 1, 
\; \nu_{pq} = \nu_{qp}$, 
and $\det (\nu_{pq}) \neq 0$. This implies that formulas (21) 
take the form 
\begin{equation}\label{eq:68}
\omega_0^a = 0,   \;\; \omega_0^p = 
g^{ps} \nu_{sq} \omega_n^q.
\end{equation}

The last equatiion in (68) show that the 1-forms 
$\omega_0^p$ are linearly independent: they are basis 
forms on the submanifold $V=V^r$ described by the point 
$A_0$ on the hyperquadric $Q^n$, on the lightlike 
hypersurface $U^n$ of rank $r$, and also on 
a frame bundle associated with this hypersurface. 
1-forms occurring in equations (4) as linear 
combinations of the basis forms $\omega_0^p$ are called 
principal forms, and the 1-forms that are not expressed 
in terms of  the basis forms are fibre forms on 
the above mentioned  frame bundle. 

By (5) the second group of equations (68) is equivalent 
to the system of equations 
\begin{equation}\label{eq:69}
\omega_p^n = \lambda_{pq}^n \omega_0^q, 
\end{equation}
where $\lambda_{pq}^n = - g_{ps} \widetilde{\nu}^{st} 
g_{ts}$,   $(\widetilde{\nu}^{st})$ is the inverse 
matrix of the matrix $(\nu_{pq})$, 
$\lambda_{pq}^n = \lambda_{qp}^n$, and $\det (\lambda_{pq}^n) 
\neq 0$. Note that we can also obtain 
equations (69) by differentiation of equation (18) 
which holds on the hypersurface $U^n$.

Taking exterior derivatives of the first 
group of equations (68), we find that 
$$
\omega_0^p \wedge  \omega_p^a = 0.
$$
Applying Cartan's lemma to this system, we find that 
\begin{equation}\label{eq:70}
\omega_p^a =  \lambda_{pq}^a  \omega_0^q, \;\; 
\lambda_{pq}^a = \lambda_{qp}^a.
\end{equation}
Note also that equations (5) and (67) imply that
$$
g_{pq} \omega_a^q + g_{ab} \omega_p^b = 0.
$$
By (70), it follows from the last equation that 
\begin{equation}\label{eq:71}
\omega^p_a =  - g_{ab} g^{pq} \lambda_{qs}^b \omega_0^s.
\end{equation}
Note also that the quantities 
$\lambda_{pq}^a$ and $\lambda_{pq}^n$ are determined in a 
second-order neighborhood of a rectilinear generator 
$l = A_0 A_n$ of the hypersurface $U^n$.

Let us prove that in our frame an $(m+1)$-dimensional span 
$L$ of the points $A_0, A_a$, and $A_n$ is a 
plane generator of the lightlike hypersurface $U^n$. 
In fact, it follows from equations (4) that 
in the case in question we have 
\begin{equation}\label{eq:72}
\renewcommand{\arraystretch}{1.3}
\left\{
\begin{array}{ll}
dA_0 =\omega_0^0 A_0      \!\!            &+ \omega_0^p A_p, \\
dA_a =\omega_a^0 A_0 + \omega_a^b A_b \!\!&+ \omega_a^p A_p 
+ \omega_a^n A_n,\\
dA_n =\omega_n^0 A_0 + \omega_n^a A_a \!\!&+ \omega_n^p A_p.
\end{array}
\right.
\renewcommand{\arraystretch}{1}
\end{equation}
If we fix the principal parameters in equations 
(72) (i.e., if we assume that $\omega_0^p = 0$), 
we obtain 
\begin{equation}\label{eq:73}
\renewcommand{\arraystretch}{1.3}
\left\{
\begin{array}{ll}
\delta A_0 =\pi_0^0 A_0, \\
\delta A_a =\pi_a^0 A_0 + \pi_a^b A_b + \pi_a^n A_n,\\
\delta A_n =\pi_n^0 A_0 + \pi_n^a A_a.
\end{array}
\right.
\renewcommand{\arraystretch}{1}
\end{equation}
In the last equations $\delta$ is the symbol of 
differentiation with respect to the fiber parameters 
(i.e., for $\omega_0^p = 0$), and $\pi_\eta^\xi 
= \omega_\eta^\xi (\delta)$.

Equations (73) show that for $\omega_0^p = 0$, 
the point $A_n$ of the hypersurface $U_n$ 
moves in an $(m+1)$-dimensional domain belonging 
to the subspace $L = A_0 \wedge  A_1 \wedge \ldots 
  \wedge  A_m \wedge  A_n$ of the same dimension. 
Let us assume that the entire subspace $L$ belongs 
to the hypersurface $U^n$, and that the point 
$A_n \in L$ moves freely in $L$. The subspace $L$ 
is tangent to the hyperquadric $Q^n$ at the point 
$A_0 \in V^r$, and thus $L$ is lightlike. 
Since the point $A_0$ describes an $r$-dimensional 
submanifold, the family of subspaces $L$ depends on 
$r$ parameters. Hence $U^n = f (M^r \times L)$, 
where $f$ is a differentiable map $f: M^r \times L 
\rightarrow P^{n+1}$. 

Equations  (72) and (73) show that the basis  
1-forms of a  lightlike  hypersurface 
$U^n$ are divided into two classes: $\omega_n^p$ and 
$\omega_n^a$. The forms  $\omega_n^p$ are connected with 
the displacement of the lightlike $(m+1)$-plane 
$L$ in the space $S^{n+1}_1$, and the forms  $\omega_n^a$ 
are connected with 
the displacement of the straight line $A_n A_0$ in this 
 $(m+1)$-plane. Since (73) implies that for $\omega_n^p = 0$ 
the point $A_0$ remains fixed, the rectilinear generator 
$A_n A_0$ describes an $m$-dimensional bundle of straight lines 
with its center at the point $A_0$, and this bundle belongs to 
the fixed $(m+1)$-dimensional subspace $L$ passing through this 
point.

Further consider an arbitrary point
\begin{equation}\label{eq:74}
z = z^0  A_0 + z^a A_a +  z^n A_n
\end{equation}
of the generator $L$ of the  lightlike  hypersurface 
$U^n$. From formulas  (72) it follows that the 
differential of any such point belongs to one and the same 
$n$-dimensional subspace $A_0 \wedge \ldots \wedge A_n$ tangent 
to the hypersurface $U^n$ at the original point $A_n$. 
The latter means that the tangent subspace to the 
hypersurface $U^n$ is not changed when the point $z$ moves along 
the lightlike generator $L$ of the hypersurface 
$U^n$. Thus,  hypersurface $U^n$ is a tangentially degenerate 
hypersurface of rank $r$. 

As a result, we arrive at the following 
theorem making Theorem 2 more 
precise:
\begin{theorem} If the rank of the tensor $\nu_{ij}$ defined by relation $(20)$ is equal to $r$, $r < n-1$, 
then a lightlike hypersurface $U^n$ of the de Sitter 
space $S^{n+1}_1$ is a ruled tangentially degenerate 
hypersurface of rank $r$ with $(m+1)$-dimensional 
lightlike generators, $m = n - r - 1$, along which  the 
tangent hyperplanes of $U^n$ are constant. The points of tangency 
of lightlike generators with the hyperquadric $Q^n$ 
form an $r$-dimensional submanifold $V^r$ on $Q^n$.
\end{theorem}

The last fact mentioned in Theorem 8 can be also treated 
in terms of quadratic hyperbands (see \cite{AG93}, p. 256). 
By Theorem 8, the hypersurface $U^n$ is the envelope of 
an $r$-parameter family of hyperplanes $\eta$ tangent 
to the hyperquadric $Q^n$ at the points of an  $r$-dimensional 
smooth submanifold $V^r$ belonging to this  hyperquadric. 
But this coincides precisely with the definition 
of quadratic hyperband. Thus Theorem 8 can be complemented as 
follows:
\begin{theorem} A lightlike hypersurface $U^n$ of rank $r$ 
in the de Sitter space $S^{n+1}_1$ is an $r$-dimensional 
quadratic hyperband with the support submanifold $V^r$ belonging 
to the Darboux hyperquadric $Q^n$.
\end{theorem}

Note also an extreme case when the rank of a  lightlike 
hypersurface $U^n$ is equal to 0. Then we have 
$$
\nu_{ij} = 0, \;\; \omega_0^i = 0.
$$
The point $A_0$ is fixed on the  hyperquadric $Q^n$, 
and the point $A_n$ moves freely in the hyperplane $\eta$ tangent 
to the  hyperquadric $Q^n$ at the point $A_0$. The 
lightlike hypersurface $U^n$ degenerates into the 
hyperplane $\eta$  tangent 
to the  hyperquadric $Q^n$ at the point $A_0$, 
and the quadratic hyperband associated with $U^n$ 
is reduced to a 0-pair consisting of the point $A_0$ and 
the hyperplane $\eta$.

Let us also find singular points on a rectilinear generator 
$L$ of a lightlike hypersurface $U^n$ of rank $r$ 
of the de Sitter 
space $S^{n+1}_1$. To this end, we write the differential of a 
point $z \in L$ defined by equation (74). We will be interested 
only in that part of this differential which does not belong to 
the generator $L$. By (72), we obtain 
$$
dz \equiv (z^0 \omega_0^p + z^a \omega_a^p + z^n \omega_n^p) A_p 
\pmod{L}.
$$
By (69), (70) and (71), we find from the last relation that 
$$
dz \equiv N^p_q (z) \omega_0^q A_p \pmod{L},
$$
where
\begin{equation}\label{eq:75}
 N^p_q (z) = \delta_q^p z^0 - g_{ab} g^{ps} \lambda^b_{sq} z^a 
- g^{ps} \lambda_{sq}^n z^n. 
\end{equation}


At singular points of a generator $L$ the dimension 
of the tangent subspace $T_x (U^n)$ to the 
hypersurface $U^n$ is reduced. By (75), this is equivalent 
to the reduction of the rank of the matrix $N^p_q (z)$. 
Thus singular points of generator $L$ can be found from the 
condition
\begin{equation}\label{eq:76}
\det N^p_q (z) = 0,
\end{equation}
which defines an algebraic focal 
submanifold ${\cal F}$ of dimension $m$ and 
order $r$ in the $(m+1)$-dimensional plane generator $L$. 
The left-hand side of equation (76) is the Jacobian of the 
map $f: M^r \times L \rightarrow P^{n+1}$ indicated above, 
and the focal submanifold ${\cal F}$ is the locus of 
singular points of this map that are located in the 
plane generator $L$ of the hypersurface $U^n$.

If the rank  of a lightlike hypersurface $U^n$ is maximal, 
that is, it is equal to $r = n - 1$, then its determinant 
manifold ${\cal F}$ is a set of singular points of its 
rectilinear generator $A_n A_0$ determined by equation (50). 
On the other hand, if $r < n -1$, then singular points of  
the straight lines $A_n A_0$ lying in the generator $L$ 
are also determined by equation (50), and they are 
the intersection points of these straight lines and the 
manifold ${\cal F}$.

In a plane generator $L$ of the lightlike hypersurface $U^n$, 
let us find an equation of the harmonic polar of 
the point $A_0$ with respect to the algebraic 
focal submanifold ${\cal F}$. Let us assume 
that the coordinates $z^a$ and $z^n$ in equation (74) 
are fixed, and the coordinate $z^0$ is variable. Then 
the point $z$ describes a straight line $l = 
A_0 \wedge (z^a A_a + z^n A_n)$. The intersection point 
of this line $l$ with the focal submanifold ${\cal F}$ 
is determined by equation (76) in which the quatities $z^a$ and 
$z^n$ are fixed, and $z^0$ is variable. Equation (76) is of 
degree $r$ with respect to $z^0$ and defines $r$ focal 
(singular) points on the straight line $l$ if each of these 
points is counted as many times as its multiplicity. 
By the Vieta theorem, the coefficient in $(z^0)^{r-1}$ in 
equation (76) is equal to the sum of roots of this equation. 
Thus, 
$$
\frac{1}{r} \sum_{p = m+1}^{n-1} z_p^0 = z^a g_{ab} \lambda^b + 
z^n \lambda^n, 
$$
where $z_p^0$ are roots of equation (76), and the quantities 
$\lambda^b$ and $\lambda^n$ are defined by the formulas 
\begin{equation}\label{eq:77}
\lambda^a = \frac{1}{r} g^{pq} \lambda^a_{pq}, \;\; 
\lambda^n = \frac{1}{r} g^{pq} \lambda^n_{pq}.
\end{equation}
Thus the harmonic pole of the point $A_0$ with respect to 
the foci of the straight line $l$ has the form
$$
C = (g_{ab} \lambda^a z^b + \lambda^n z^n) A_0 + z^a A_a + z^n A_n.
$$
The locus of these poles on a 
plane generator $L$ of the hypersurface $U^n$ is defined by the 
equation
\begin{equation}\label{eq:78}
z^0 -  g_{ab} \lambda^a z^b - \lambda^n z^n = 0, 
\end{equation}
whose left-hand side is the trace of the matrix $N_q^p (z)$. 
Equation (78) defines a subspace of dimension $m$ on 
the $(m+1)$-dimensional generator $l$ of 
the hypersurface $U^n$. This subspace is the harmonic polar 
of the point $A_0$ with respect to the algebraic 
focal submanifold ${\cal F}$.

For construction of screen distribution on a lightlike 
hypersurface $U^n$ of rank $r < n - 1$ we will need differential 
prolongations of Pfaffian equations (69) and (70). 
Taking exterior derivatives of these equations, we find that 
\begin{equation}\label{eq:79}
(\nabla \lambda^a_{pq} + \lambda^a_{pq} \omega_0^0 
 + \lambda^n_{pq} \omega_n^0 + g_{pq} \omega_{n+1}^a) \wedge 
 \omega_0^q = 0,  
\end{equation}
\begin{equation}\label{eq:80}
(\nabla \lambda^n_{pq} + \lambda^n_{pq} \omega_0^0 
 + \lambda^a_{pq} \omega_a^n + g_{pq} \omega_{n+1}^n) \wedge 
 \omega_0^q = 0,  
\end{equation}
where 
$$
\renewcommand{\arraystretch}{1.3}
\begin{array}{ll}
\nabla \lambda^a_{pq} = d \lambda^a_{pq} 
- \lambda^a_{sq} \omega_p^s - \lambda^a_{ps} \omega_q^s 
+ \lambda^b_{pq} \omega_b^a, \\
\nabla \lambda^n_{pq} = d \lambda^n_{pq} 
- \lambda^n_{sq} \omega_p^s - \lambda^n_{ps} \omega_q^s.
\end{array}
\renewcommand{\arraystretch}{1}
$$
Applying Cartan's lemma to equations (79) and (80) 
and fixing the principal parameters (i.e., setting 
$\omega_0^p = 0$), we find that 
\begin{equation}\label{eq:81}
\renewcommand{\arraystretch}{1.3}
\left\{
\begin{array}{ll}
\nabla_\delta \lambda^a_{pq} + \lambda^a_{pq} \pi_0^0 
+ \lambda^a_{pq} \pi_n^0 + g_{pq} \pi_{n+1}^a = 0, \\
\nabla_\delta \lambda^n_{pq} + \lambda^n_{pq} \pi_0^0 
+ \lambda^a_{pq} \pi_a^n + g_{pq} \pi_{n+1}^n = 0.
\end{array}
\right.
\renewcommand{\arraystretch}{1}
\end{equation}
Note also that by the last equation of equations (5),
 the tensor $g_{pq}$ defined by the first group of relations 
(67) satisfies the equations
\begin{equation}\label{eq:82}
\nabla_\delta g_{pq}= 0.  
\end{equation}
Equations (81) and (82) prove that neither quantities 
$\lambda^a_{pq}$ nor quantities $\lambda^n_{pq}$ form 
a geometric object, but jointly, the quantities  
$\lambda^a_{pq}, \lambda^n_{pq}$, and 
the tensor $g_{pq}$ form a linear geometric object. 

By equations (81) and (82), the quantities $\lambda^a$ and 
$\lambda^n$ defined by formulas (77) satisfy the equations
\begin{equation}\label{eq:83}
\renewcommand{\arraystretch}{1.3}
\left\{
\begin{array}{ll}
\nabla_\delta \lambda^a + \lambda^a \pi_0^0 
+ \lambda^n \pi_n^a + \pi_{n+1}^a = 0, \\
\nabla_\delta \lambda^n +  \lambda^n \pi_0^0 
+ \lambda^a \pi_a^n + \pi_{n+1}^n = 0.
\end{array}
\right.
\renewcommand{\arraystretch}{1}
\end{equation}
It follows that jointly, the quantities $\lambda^a$ and 
$\lambda^n$ form a geometric object which is associated with 
a second-order differential neighborhood of a plane generator $L$ 
of the hypersurface $U^n$.

Next we construct the quantities
\begin{equation}\label{eq:84}
a_{pq}^a = \lambda_{pq}^a - \lambda^a g_{pq}, \;\; 
a_{pq}^n = \lambda_{pq}^n - \lambda^n g_{pq}.  
\end{equation}
By (80),(81), and (82), they satisfy the equations
\begin{equation}\label{eq:85}
\renewcommand{\arraystretch}{1.3}
\left\{
\begin{array}{ll}
\nabla_\delta a_{pq}^a + a_{pq}^a \pi_0^0 + a_{pq}^n \pi_n^0 
= 0, \\ 
\nabla_\delta a_{pq}^n + a_{pq}^n \pi_0^0 + a_{pq}^a \pi_a^n 
= 0. 
\end{array}
\right.
\renewcommand{\arraystretch}{1}
\end{equation}
These equations prove that jointly the quantities $a_{pq}^a$ 
and $a_{pq}^n$ form a tensor with respect to 
the admissible transformations in a frame bundle associated 
with the hypersurface $U^n$. Let us assume that the indices 
$\alpha$ and $\beta$ take $m + 1$ values $1, \ldots , m, n$. 
Then the equations which the tensor $a_{pq}^\alpha = \{a_{pq}^a, 
a_{pq}^n\}$ satisfy can be written in the form 
\begin{equation}\label{eq:86}
\nabla_\delta a_{pq}^\alpha + a_{pq}^\alpha \pi_0^0 = 0,
\end{equation}
where $\nabla_\delta a_{pq}^\alpha = \delta a_{pq}^\alpha 
- a_{sq}^\alpha \pi^s_p - a_{ps}^\alpha \pi^s_q 
+ a_{pq}^\beta \pi_\beta^\alpha$. 
It follows from equations (84) that this tensor satisfies 
the apolarity condition 
\begin{equation}\label{eq:87}
a_{pq}^\alpha g^{pq} = 0.
\end{equation}

Consider the straight lines $A_0 A_\alpha$ connecting 
the point $A_0$ with the points $A_\alpha, \; \alpha = 1, 
\ldots , m, n$. Their parametric equations can be written in the 
form 
$$
\widetilde{A}_\alpha = A_\alpha + x_\alpha A_0.
$$
Let us find the points of intersection of 
these straight lines with 
the harmonic polar $K$ of the point $A_0$ 
with respect to the focal submanifold ${\cal F}$. 
Substituting coordinates of these points into equation (78), we 
find that
$$
x_\alpha = \lambda_\alpha, \;\;\;\; \alpha = 1, 
\ldots , m, n,
$$
where $\lambda_a = g_{ab} \lambda^b$ and $\lambda_n = \lambda^n$. 

The points $C_\alpha = A_\alpha + \lambda_\alpha A_0$ lying in 
the subspace $K$ can be taken as the vertices of a reduced 
frame associated with the hypersurface $U^n$ and defined in 
a second-order differential neighborhood of the plane generator 
$L$ of this hypersurface. If we consider our hypersurface 
with respect to this reduced frame, then we have 
\begin{equation}\label{eq:88}
\lambda^\alpha = 0, \;\; \lambda^\alpha_{pq} = a^\alpha_{pq}.
\end{equation}
It follows from equation (83) that the 1-forms 
$\omega_{n+1}^\alpha$ become principal forms:
\begin{equation}\label{eq:89}
\omega_{n+1}^\alpha = b_p^\alpha \omega^p.
\end{equation}
With respect to  new frame, equations (69) and (70) take the form 
\begin{equation}\label{eq:90}
\omega_p^\alpha = a_{pq}^\alpha \omega^q, \;\; a_{pq}^\alpha = 
a_{qp}^\alpha.
\end{equation}

Consider the rectangular matrix $A = (a_{pq}^\alpha)$ in which 
 $\alpha$ is the row number, and the pair $(p, q) = (q, p)$ is 
the column number. The matrix $A$ has $m + 1$ rows and 
$\frac{1}{2} r (r + 1)$ columns. But by (87), not more than 
$\frac{1}{2} r (r + 1) - 1$ columns of the matrix $A$ are 
linearly independent. Suppose that $\mbox{{\rm rank}} A = \rho, 
\; \rho \leq \mbox{{\rm min}} \{m+1, \frac{1}{2} r (r + 1) - 1\}$. 
Construct the following tensors:
\begin{equation}\label{eq:91}
a^{\alpha\beta} = g^{pq} g^{st} a^\alpha_{ps} a^\beta_{qt} 
\;\;\;\mbox{{\rm and}} \;\;\; a_\beta^\alpha = 
g_{\beta\gamma} a^{\gamma\alpha}.
\end{equation}
It is not difficult to prove that the rank of each of 
these tensors is equal to 
the rank of the matrix $A$, $\mbox{{\rm rank}} (a^{\alpha\beta}) 
= \mbox{{\rm rank}} (a_\beta^\alpha) = \rho$.

Construct the quantity
$$
a = a_{[\alpha_1}^{\alpha_1} a_{\alpha_2}^{\alpha_2} \ldots 
a_{\alpha_\rho]}^{\alpha_\rho},
$$
which is equal to the sum of the diagonal minors of order $\rho$ 
of the matrix $(a^\alpha_\beta)$. Since the rank of this matrix 
is equal to $\rho$, then if $\rho \geq 1$, the quantity $a$ is 
different from 0, $a \neq 0$.

From equations (86) it follows that the tensor $a_\beta^\alpha$ 
satisfies the equations
$$
\nabla_\delta a_\beta^\alpha + 2 a_\beta^\alpha \pi_0^0 = 0.
$$
Applying the formula for differentiation of determinants, 
we find from the last equation that the quantity    $a$ 
satisfies the equation
\begin{equation}\label{eq:92}
\delta a + 2 \rho a \pi_0^0 = 0,
\end{equation}
i.e., $a$ is a relative invariant of weight $- 2\rho$.

Equation (92) is written for fixed principal parameters, i.e., 
under the condition $\omega_0^p = 0$. If these parameters 
are variable, then it follows from equation (92) that 
\begin{equation}\label{eq:93}
\frac{da}{2\rho a} +  \omega_0^0 = \mu_p \omega_0^p.
\end{equation}
Taking the exterior derivative of the last equation, we find that 
\begin{equation}\label{eq:94}
(d\mu_p - \mu_q \omega^q_p + \omega_p \omega_0^0 + \omega_p^0) 
 \wedge \omega_0^p = 0.
\end{equation}
This implies that the quantities $\mu_p$ form a geometric 
object. For $\omega_0^p = 0$, this object satisfies the 
equations
\begin{equation}\label{eq:95}
\nabla_\delta \mu_p + \mu_p \pi_0^0 + \pi_p^0 = 0.
\end{equation}
It follows from equation (93) that the geometric object $\mu_p$ 
is defined in a third-order differential neighborhood of the 
plane generator $L$ of the hypersurface $U^n$.

Consider the subspace $T_{A_0} (V^r) 
= A_0 \wedge A_{m+1} \wedge \ldots \wedge A_{n-1}$ 
tangent to the submanifold $V^r$ described by 
the point $A_0$ on the hyperquadric $Q^n$. This subspace 
belongs to the tangent hyperplane $\eta$ to the 
lightlike hypersurface $U^n$ and is orthogonal to its 
plane generator $L =  A_0 \wedge A_1 \wedge \ldots \wedge A_m 
\wedge A_n$. The geometric object $\mu_p$ allows us to 
construct a normalizing subspace $\zeta$ 
of dimension $r-1$ in  $T_{A_0} (V^r)$. To this end, consider 
the points
$$
\widetilde{A}_p = A_p + x_p A_0.
$$
Differentiating these points, applying equations 
(5), (69), and (70), and assuming that the 
principal parameters are fixed, i.e., 
$\omega^p_0 = 0$, we find that 
$$
\delta \widetilde{A}_p = (\nabla_\delta x_p + x_p \pi_0^0 + \pi_p^0) A_0 + \pi_p^q \widetilde{A}_q.
$$
This implies that the subspace spanned by the points 
$\widetilde{A}_p$ is invariant if and only if the quantities 
$x_p$ satisfy the differential equations
$$
\nabla_\delta x_p + x_p \pi_0^0 + \pi_p^0 = 0.
$$
Comparing these equations with equations (95), we see that 
they are satisfied if we take $x_p = \mu_p$. Thus the points
$$
C_p = A_p + \mu_p A_0
$$
determine an invariant normalizing subspace $\zeta = 
C_{m+1} \wedge \ldots \wedge C_{n-1}$.

Suppose that $x$ is an arbitrary point of the generator 
$L$ of the hypersurface $U^n$. This point and the subspace 
$\zeta$ define an $r$-plane $\Delta_x = x \wedge \zeta$. 
Such $r$-planes are defined for all points $x \in U^n$ and 
form an $r$-dimensional {\em screen distribution} $\Delta$ on 
$U^n$ which is complementary to the generators $L$ of $U^n$. 
Since the geometric object $\mu_p$ is defined in a third-order 
neighborhood, then the screen distribution is defined in the same 
neighborhood. Thus the following theorem holds.

\begin{theorem}
 If the rank of the matrix $A = (a^\alpha_{pq})$ is different 
from $0$, then in a third-order 
neighborhood of a plane generator $L$ of 
a lightlike hypersurface $U^n$ of rank 
$r < n - 1$, there is defined an invariant $r$-dimensional screen 
distribution $\Delta$.
\end{theorem}

If we place the vertices $A_p$ of our frame into the 
points $C_p$, then we obtain $\mu_p = 0$. This and equation 
(94) imply that 
\begin{equation}\label{eq:96}
 \omega_p^0 = c_{pq}  \omega_0^q, \;\; c_{pq} = c_{qp}.
\end{equation}
The points $C_p$ and the normalizing subspace 
$\zeta = C_{m+1} \wedge \ldots \wedge C_{n-1}$ 
are defined in  a third-order 
neighborhood  of a plane generator $L$, and the quantities 
$c_{pq}$ are defined in a fourth-order neighborhood.

Let us prove that the  
fibration of normalizing subspaces we have constructed
determines an affine connection on a hypersurface $U^n$ 
which can be  considered as an $r$-parameter fibration of its 
$(m+1)$-dimensional plane generators $L$. In fact, the basic forms of this fibration are the 1-forms $\omega_0^p$. 
Taking exterior derivatives of these forms, we find that
\begin{equation}\label{eq:97}
d \omega^p_0 =  \omega_0^q \wedge \theta_q^p, 
\end{equation}
where $\theta_q^p = \omega_q^p - \delta_q^p \omega_0^0$. 
Taking exterior derivatives of the forms $\theta_q^p$ and taking into account that by (96) $d\omega_0^0 = 0$, we obtain
\begin{equation}\label{eq:98}
d  \theta_q^p -  \theta_q^s \wedge \theta_s^p  
= R^p_{qst}  \omega_0^s \wedge  \omega_0^t, 
\end{equation}
where 
\begin{equation}\label{eq:99}
R^p_{qst} = - g_{\alpha\beta} g^{pu} a^\alpha_{q[s} a^\beta_{t]u} 
+ c_{q[s} \delta^p_{t]} + g_{q[s} g^{pu} c_{t]u},
\end{equation}
and
$$
g_{\alpha\beta} = \pmatrix{g_{ab} & 0 \cr
                               0  & 1 \cr}.
$$

Thus the following theorem is valid.

\begin{theorem}
An invariant screen distribution induces a torsion-free 
affine connection on the fibration of plane generators of a 
hypersurface $U^n$. The curvature tensor of this connection is 
determined by equation $(99)$, and its tensor Ricci is symmetric.
\end{theorem}

The last statement of Theorem 11 can be proved by a direct 
calculation. In fact, contracting equation (99) with respect 
to the indices $p$ and $t$, we find that 
$$
R_{pq} = R^s_{pqs} = \frac{1}{2} \Bigl(g_{\alpha\beta} g^{st} 
a^\alpha_{ps} a^\beta_{qt} + (r-2) c_{pq} 
+ g_{pq} g^{st} c_{st}\Bigr),
$$
It is easy to see that this tensor is symmetric, 
$R_{pq} = R_{qp}$.

We did not consider yet only the case when the rank of the matrix 
$A = (a_{pq}^\alpha)$ is equal to 0, $\rho = 0$. In this case 
the matrix $A$ is the null-matrix, and the construction 
of an invariant screen distribution is impossible. Let us clarify 
the geometric structure of the hypersurface $U^n$ in this case.

If this is the case, formulas (84) imply that 
\begin{equation}\label{eq:100}
\lambda^a_{pq} = \lambda^a g_{pq}, \;\;
\lambda_{pq}^n = \lambda^n  g_{pq}.
\end{equation}
Thus the Jacobi matrix (75) of 
 the mapping $f: M^r \times L \rightarrow P^{n+1}$ takes 
the form
$$
N_q^p (z) = \delta_q^p (z^0 - g_{ab} \lambda^b - \lambda^n),
$$
and the equation of the focal submanifold ${\cal F}$ becomes
$$
\det N_q^p (z) = (z^0 - g_{ab} \lambda^b z^a - z^n)^r = 0.
$$
Thus the focal submanifold ${\cal F}$  is an $r$-fold linear 
subspace
\begin{equation}\label{eq:101}
z^0 - g_{ab} \lambda^b z^a - z^n = 0
\end{equation}
of dimension $m$ belonging to the $(m+1)$-dimensional plane 
generator $L$ of the hypersurface $U^n$. It is possible to prove 
that if $r \geq 2$, then this subspace is fixed, and 
the hypersurface $U^n$ is an $n$-dimensional cone with 
an $(m+1)$-dimensional plane generators and an $m$-dimensional 
vertex defined by equation (101).

In this case the submanifold $V^r$, along which the hypersurface 
$U^n$ is tangent to the hyperquadric $Q^n$, is an $r$-dimensional
 sphere $S^r$.

If $r = 1$, then the hypersurface 
$U^n$ is an envelope of a one-parameter family of isotropic 
hyperplanes that are tangent to the hyperquadric $Q^n$ at 
the points of an arbitrary curve $\gamma$. 
Finally if $r = 0$, then the hypersurface 
$U^n$ is an isotropic hyperplane.

Note that the invariant normalization of a 
lightlike hypersurface $U^n$ which we have constructed in 
Section  6 is a new geometric interpretation of a 
conformally invariant normalization of a submanifold $V^r$ 
of a conformal space $C^n$ which was constructed in 
\cite{A61} (see also \cite{AG96}, Ch. 3).

{\sf             
Department of Mathematics,            
Jerusalem College of Technology - Mahon Lev, 
21 Havaad Haleumi St., 
P.O. Box 16031,  Jerusalem 91160, Israel       
}

\noindent
{\em E-mail address}: akivis@math.jct.ac.il

\vspace*{3mm}

\noindent
{\sf 
 Department of Mathematics,   New Jersey Institute of Technology,  
 University Heights, Newark, NJ 07102
}

\noindent
{\em E-mail address}: vlgold@numerics.njit.edu


\begin{thebibliography}{AG96}


\bibitem{A52} Akivis, M. A.,
 {\em Invariant construction of the geometry of a 
  hypersurface of a conformal space},   Mat. Sb. (N.S.)
  {\bf 31} (73) (1952), no. 1, 43--75 (Russian). 



\bibitem{A61} Akivis, M. A.,
 {\em On the conformal differential geometry of 
  multidimensional surfaces},  Mat. Sb. (N.S.) {\bf 53}
 (95) (1961), no. 4, 53--72 (Russian). 

\bibitem{A63} Akivis, M. A.,
  {\em On the structure of multidimensional surfaces carrying
  a net of curvature lines},  Dokl. Akad. Nauk SSSR
 {\bf 149} (1963), no. 6, 1247--1249 (Russian); 
English transl. in Soviet 
 Math. Dokl. {\bf 4} (1963), no. 2, 529--531. 


\bibitem{A64} Akivis, M. A., 
{\em Conformal differential geometry},  Geometry
  1963, pp. 108--137. Akad. Nauk SSSR Inst. Nauchn. Informatsii,
      Moscow, 1965 (Russian). 


\bibitem{A65} Akivis, M. A.,
 {\em On an invariant differential-geometric characterization 
 of the Dupin cyclide},  Uspekhi Mat. Nauk 
 {\bf 20} (1965), no. 1, 177--180 (Russian). 



\bibitem{AG93}  Akivis, M. A. and V. V. Goldberg,
 {\em Projective differential geometry of submanifolds.} 
North-Holland, Amsterdam-New York-Tokyo, 1993, xi+364 pp. 


\bibitem{AG96}  Akivis, M. A. and V. V. Goldberg,
 {\em Conformal differential geometry and its generalizations}. 
 Wiley-Interscience Publication, New York,   1996, xiv+383 pp.


\bibitem{AG97}  Akivis, M. A. and V. V. Goldberg,
 {\em Darboux mapping of canal hypersurfaces} 
(submitted),   1997.



\bibitem{C50} Casanova, G:
{\em La notion de p\^{o}le harmonique}, Rev. math. spec. 
{\bf 65} (1950), no. 6, 437--440.


\bibitem{DB96} Duggal, K. L., and A. Bejancu, 
{\em Lightlike submanifolds of semi-Riemannian manifolds 
and applications}, Kluwer Academic Publishers, Amsterdam, 
 1996, 308 pp.




\bibitem{F56} Finikov, S. P.,  
{\em Theory of pairs of congruences}, Gosudarstv. Izdat. 
Tehn.-Teor. Lit., Moscow, 1956, 443 pp. (Russian);  French transl.  by M. Decuyper, U.E.R. 
Mathematiques Pures et Appliquees,  No. 68,   Universit\'{e} des
 Sciences et Techniques de Lille I, Villeneuve d'Ascq,  
$n^\circ 68$, 2 vols, 1976, xxix+616 pp.  

\bibitem{K89} Kossowski, M., 
{\em The intrinsic conformal structure and Gauss map of 
a light-like hypersurface in Minkowski space}, 
  Trans. Amer. Math. Soc.  {\bf 316} (1989), no. 1,  369--383.   


\bibitem{ON83} O'Neill,  B.,  
{\em Semi-Riemannian geometry. With applications to relativity}, Academic Press, New York, 1983, 
xiii+468 pp.  


\bibitem{W77} Wolf, J. A.,
{\em Spaces of constant curvature}, 4th ed., Publish or Perish, 
Berkeley,  1977, xvi+408 pp. 


\bibitem{Z96} Zheng, Y., 
{\em Space-like hypersurfaces with constant scalar curvature in the de Sitter space}, 
  Differential Geom.  Appl. {\bf 6} (1996),  51--54.   



\end{thebibliography}
\end{document}